%% file: main.tex
\documentclass[preprint,12pt,authoryear]{elsarticle}
\pdfoutput=1

\usepackage{amssymb}

\RequirePackage{fix-cm}

\usepackage{parskip}

\usepackage{graphicx}
\usepackage[justification=centering]{caption}

\usepackage{mathptmx}

\usepackage[utf8]{inputenc}

\usepackage{amsmath,empheq}
\usepackage{amssymb}
\usepackage[english]{babel}
\usepackage{epsfig}
\usepackage{multirow,multicol}
\usepackage{algorithmic}
\usepackage{algorithm}
\usepackage{mathrsfs}
\usepackage{dsfont}
\usepackage{appendix}
\usepackage{comment}

\usepackage{tikz}
\usepackage{tikz-cd}

\usepackage{
,nicefrac 
,pgf
,geometry 
,graphics
,graphicx 
,pdfpages 
,pgfplots
}

\usetikzlibrary{decorations.pathreplacing,quotes}
\usetikzlibrary{shapes}
\usetikzlibrary{positioning,chains,fit,calc}
\usetikzlibrary{backgrounds}

\tikzstyle{vertex}=[scale=1,auto=left,align=center, style={circle,fill=blue!20}]
\tikzstyle{blank}=[scale=1,auto=left,align=center, style={circle}]
\tikzstyle{comment}=[scale=1,auto=left,align=center, style={circle,fill=white!20}]

\tikzstyle{sommet}=[ draw, circle, fill=blue!20, minimum size=6mm, inner sep=0pt]
\tikzstyle{petitsommet}=[ draw, circle, fill=blue!20, minimum size=3mm, inner sep=0pt]
\tikzstyle{sommetFeu}=[ draw, star,star points=10, fill=red!20, minimum size=6mm,inner sep=0pt]
\tikzstyle{centre}=[ draw, regular polygon,regular polygon sides=5, fill=blue!20, minimum size=6mm, inner sep=0pt]
\tikzstyle{smallsepcentre}=[ draw, regular polygon,regular polygon sides=5, fill=blue!20, minimum size=6mm,inner sep=0pt]
\tikzstyle{centreFeu}=[ draw, regular polygon,regular polygon sides=5, fill=red!20, minimum size=6mm,inner sep=0pt]

\tikzstyle{sommetpt}=[ draw, circle, fill=blue!20, minimum size=3mm, inner sep=0pt]

\usepackage{fontawesome}

\usepackage[T1]{fontenc}

\usepackage{multirow}
\usepackage{booktabs}

\usepackage{array}
\usepackage{tabularx}
\usepackage{etoolbox}
	\appto\TPTnoteSettings{\footnotesize}

\newtheorem{theoreme}{Theorem}
\newtheorem{corollary}{Corollary}

\newtheorem{remarque}{Remark}
\newtheorem{monclaim}{Claim}
\newtheorem{transformation}{Transformation}
\newtheorem{example}{Example}
\newtheorem{lem}{Lemma}
\newproof{pf}{Proof}
\newtheorem{proposition}{Proposition}

\newcommand{\bgenum}{\begin{enumerate}}
\newcommand{\edenum}{\end{enumerate}}

\newcommand{\bgitem}{\begin{itemize}}
\newcommand{\editem}{\end{itemize}}

\DeclareMathOperator*{\argmax}{\arg\!\max\ }
\DeclareMathOperator*{\argmin}{\arg\!\min\ }

\usepackage{environ}

\makeatletter
\newcommand{\problemtitle}[1]{\gdef\@problemtitle{#1}}
\newcommand{\probleminput}[1]{\gdef\@probleminput{#1}}
\newcommand{\problemfeasible}[1]{\gdef\@problemfeasible{#1}}
\newcommand{\problemoptimization}[1]{\gdef\@problemoptimization{#1}}
\NewEnviron{problem}{
  \problemtitle{}\probleminput{}\problemfeasible{}\problemoptimization{}
  \BODY
  \par\addvspace{.5\baselineskip}
  \noindent
  \begin{tabularx}{\textwidth}{@{\hspace{\parindent}} l X c}
    \multicolumn{2}{@{\hspace{\parindent}}l}{\@problemtitle} \\
    \textbf{Instance:} & \@probleminput \\
    \textbf{Feasible solutions:} & \@problemfeasible\\
    \textbf{Objective:} & \@problemoptimization
  \end{tabularx}
  \par\addvspace{.5\baselineskip}
}
\makeatother

\usepackage{todonotes}

\newcommand\mar[1]{\textcolor{black}{#1}}
\newcommand\had[1]{\textcolor{black}{#1}}
\newcommand\cec[1]{\textcolor{black}{#1}}
\pgfplotsset{compat=1.16}

\journal{arxiv.org}

\begin{document}

\begin{frontmatter}



\title{Hardness and approximation of the {Probabilistic $p$-Center} problem {under Pressure}\tnoteref{t1}}
 \tnotetext[t1]{Work supported by the European Union’s Horizon 2020 research and innovation programme under the Marie Skłodowska-Curie grant agreement No 691161.}
 

\author[1]{Marc Demange\corref{cor1}}
\ead{marc.demange@rmit.edu.au}

\author[1,2]{Marcel A. Haddad}
\ead{marcel.haddad@dauphine.fr}


\author[2]{Cécile Murat}
\ead{cecile.murat@dauphine.fr}

 \cortext[cor1]{Corresponding author}

 \address[1]{School of Science, RMIT University, Melbourne, Vic., Australia}
 \address[2]{Universit\'e Paris-Dauphine, PSL Research University, CNRS, LAMSADE, 75016 Paris, France}

\begin{abstract}
The Probabilistic $p$-Center problem under Pressure ({\tt Min P$p$CP}) is a variant of the usual {\tt Min $p$-Center} problem we recently introduced in the context of wildfire management. The problem is 
to
locate $p$ shelters minimizing the maximum distance people will have to cover 
to reach the closest accessible shelter in case of fire. 
The landscape is divided \had{into} 
zones and is modeled as an edge-weighted graph with vertices corresponding to zones and edges corresponding to direct connections between two adjacent zones.
The risk associated with fire outbreaks is modeled using a finite set of fire scenarios. Each scenario 
corresponds to a fire outbreak on a single zone (i.e., on a vertex) with the main consequence of modifying evacuation paths in two ways. First, an evacuation path cannot pass through the vertex on fire.
Second, the fact that 
someone close to the fire may not take rational decisions when selecting a direction to escape is modeled using  new kinds of evacuation paths. In this paper, for a given instance of {\tt Min P$p$CP} defined by an edge-weighted graph $G=(V,E,L)$ and an integer $p$, we characterize the set of feasible solutions of {\tt Min P$p$CP}. We prove that {\tt Min P$p$CP} cannot be approximated with a ratio less than $\frac{56}{55}$ on  subgrids (subgraphs of grids) of degree at most 3. Then, we propose some approximation results for {\tt Min P$p$CP}. These results require approximation results for two variants of the (deterministic) {\tt Min $p$-Center} problem called {\tt Min MAC $p$-Center} and {\tt Min Partial $p$-Center}. 
\end{abstract}




\end{frontmatter}





\section{Introduction}

The problem {\tt Min P$p$CP} was introduced in~\cite{DHM18} as a variant of the usual {\tt Min $p$-Center} problem with 
indeterminacy on vertices. In the same paper, we presented  our motivation in the context of wildfire management and discussed it further in~\cite{demange2020robust}.
In our model, the landscape is represented by an adjacency graph $G=(V,E)$. Each vertex corresponds to a zone and two vertices $i$ and $j$ are linked by an edge if and only if it is possible to go directly from one to the other without passing through another area. We assume this is a symmetric relation, which makes this graph non-directed. Each edge $(i,j)$ is weighted
with a positive number $\ell_{ij}$ that can be seen as a distance or a traveling time; we will call it {\em the length} of the edge $(i,j)$.  
For every two vertices $i,j$, $d(i,j)$ will denote the shortest path distance between $i$ and $j$ in $G$ and for any set of vertices $C\subset V$, we denote $d(v,C)=\min\limits_{c\in C} d(v,c)$ the distance from $v$ to $C$. By convention, we will set $d(v,\emptyset)=+\infty$.

For a given integer $p$, the objective is to select a set $C$ (called $p$-center) of at most $p$ vertices, i.e. zones, where to locate fire-proof shelters so as to minimize the maximum traveling time from a zone to a shelter. In a deterministic setup, this problem is the classical {\tt Min $p$-Center} problem that aims to locate facilities on vertices of a network modeled by a graph. Given our motivating context,  centers will just be called shelters and, when no ambiguity occurs, we will just use the term shelter to refer to a vertex where to install a shelter. For a set $C$  of shelters and a vertex $j$, $d(j,C)$ will be called {\em {distance to shelters}} of  $j$ and the {(deterministic}) {\em radius} of $C$, denoted $r(C)$,  corresponds to the longest {distance to shelters} of vertices: $r(C)=\max\limits_{v\in V}d(v,C)$. {\tt Min $p$-Center} is to find, for any $p$,  a set $C, |C|\leq p$ of minimum radius. 

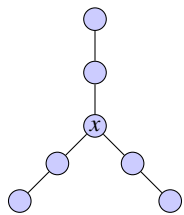
\begin{figure}[hbt]
\centering
\begin{tikzpicture}

\node[petitsommet] (x) at (0,0){$\scriptstyle x$};
\node[petitsommet] (a) at (0,0.71){};
\node[petitsommet] (b) at (0,1.42){};
\node[petitsommet] (c) at (0.5,-0.5){};
\node[petitsommet] (d) at (1,-1){};
\node[petitsommet] (e) at (-0.5,-0.5){};
\node[petitsommet] (f) at (-1,-1){};

\foreach \from/\to in {x/a,a/b,x/c,c/d,x/e,e/f}
 \draw (\from) -- (\to); 

\end{tikzpicture}
\caption{An example where, for $p=2$, a singleton minimizes the radius.}
\label{fig:star}
\end{figure}

Since adding a center to $C$ cannot increase its radius, it is straightforward that, if  $p\leq |C|$, then there is an optimal solution with exactly $p$ shelters; however, this is not a necessary condition for optimal solutions. Consider indeed the graph of Figure~\ref{fig:star} with all edge lengths equal to~1; if $p=2$, then the minimum possible radius is~2 but $r(\{x\})=2$. 

{\tt Min $p$-Center} and numerous versions   have been extensively studied both from a graph theory perspective and for various applications~(see, for instance~\cite{calik2015p}). 
It is a well known NP-hard problem, even in the class of planar graphs \had{with degrees less than~3} (\cite{KH79}) that is particularly relevant in our motivating context.  {\tt Min $p$-Center} is known to be 2-approximable~(\cite{HS85}) and is not approximated with a constant ratio strictly smaller than~2, unless P=NP~(\cite{Hard-Bottleneck}). Similar results can be obtained for variants of {\tt Min $p$-Center}. For instance, 
in \cite{chaudhuri1998p},  the 
generalization of {\tt Min $p$-Center} where, given a number $k$, we have to place $p$ centers so as to minimize the maximum distance of any non-center node to its $k^{th}$ closest center. They give a 2-approximation algorithm for this problem, and show it is the best possible. In this paper, \mar{to establish approximation results for the problem we  deal with,} we will 
need to establish similar results for two variants of {\tt Min $p$-Center}: {\tt Min MAC $p$-Center} defined in Section~\ref{sec:approximation} and {\tt Min Partial $p$-Center} introduced in~\cite{Partial-p-center}.

{\tt Min P$p$CP},  is a version of {\tt Min $p$-Center} with indeterminacy on vertices: with some probability, a vertex may become unavailable due to a fire outbreak. We present this problem in details in Section~\ref{sec:pb} after giving required related definitions. We discuss the difference with other versions of {\tt Min $p$-Center} under indeterminacy and characterize feasible solutions. 
Then, in Section~\ref{sec:complexity}, we present a new hardness result in instance classes that are natural for our motivating application.
Finally, in Section~\ref{sec:approximation}, we investigate some approximation results. To this purpose we will use {\tt Min MAC $p$-Center} and {\tt Min Partial $p$-Center}. The former problem aims at finding a  $p$-center that is a feasible solution for the problem {\tt Min P$p$CP} (see Section~\ref{ssec:feasible_solutions}) and of minimum radius.
The latter problem is to find a $p$-center of minimum partial radius, where only some vertices are taken into account to compute the partial radius (only these vertices are required to be close to a center).

All definitions of problems used in the paper are recalled in~\ref{app:listproblems}. Main hardness and approximation results are reported in Table~\ref{table-result}. 

\begin{table}[h]
\centering
\footnotesize
\begin{tabular}{|c|c|c|c|}
\hline 

\multirow{2}{*}{\bf Instance class}  & \multirow{2}{*}{\bf Complexity} &\multicolumn{2}{c|}{\bf Approximation} \\
  &  &{Lower bound} & {Upper bound} \\
  \hline
  \multicolumn{4}{|c|}{\multirow{2}{*}{\tt Min P$p$CP}}\\
  \multicolumn{4}{|c|}{}\\
 \hline
\multirow{2}{*}{Bipartite planar } & \multicolumn{1}{c|}{\multirow{2}{*}{NP-hard}} & \multirow{2}{*}{$\frac{20}{19}$}& \multirow{2}{*}{\bf 14}\\
&\multicolumn{2}{c|}{}&\\
of degree~2 or 3&\multicolumn{2}{c|}{\em \cite{DHM18}}&{\em Th.~\ref{cor:approx}}\\
\hline 
\multirow{2}{*}{Subgrid of deg. at most 3} & \multicolumn{1}{c|}{\multirow{2}{*}{\bf{NP-hard}}} & \multirow{2}{*}{$\mathbf{\frac{56}{55}}$}& \multirow{2}{*}{\bf 14}\\
&\multicolumn{2}{c|}{}&\\
and all edge lengths 1 &\multicolumn{2}{c|}{\em Th.~\ref{thm:complexity}}&{\em Th.~\ref{cor:approx}}\\
\hline 
\multirow{2}{*}{Tree} & \multicolumn{2}{c|}{\multirow{2}{*}{?}}&\multirow{2}{*}{\bf 3}\\
&\multicolumn{2}{c|}{}&\\
and all edge lengths 1 &\multicolumn{2}{c|}{ }&{\em Cor.~\ref{cor:tree}}\\
\hline 
\multirow{2}{*}{Any graph $G$ with} & \multicolumn{1}{c|}{\multirow{2}{*}{NP-hard}} & \multirow{2}{*}{$\frac{20}{19}$}& \multirow{2}{*}{$\mathbf{4\overline{deg}(G) +2}$}\\
&\multicolumn{2}{c|}{}&\\
lengths in $[\ell, 2\ell]$&\multicolumn{2}{c|}{\em \cite{DHM18}}&{\em Th.~\ref{cor:approx}}\\
\hline
 \multicolumn{4}{|c|}{\multirow{2}{*}{\tt Min MAC $p$-Center}}\\
  \multicolumn{4}{|c|}{}\\
 \hline
\multirow{2}{*}{Any} & \multicolumn{1}{c|}{\multirow{2}{*}{\bf{NP-hard}}} & \multirow{2}{*}{$\mathbf{2}$}& \multirow{2}{*}{$\mathbf{2}$}\\
&\multicolumn{2}{c|}{}&\\
&\multicolumn{2}{c|}{\em Prop.~\ref{prop:lb2mac}}&{\em Th.~\ref{th: approx-mac_pcenter}}\\
 \hline
   \multicolumn{4}{|c|}{\multirow{2}{*}{\tt Min Partial $p$-Center}}\\
  \multicolumn{4}{|c|}{}\\
  \hline
\multirow{2}{*}{Any} & \multicolumn{1}{c|}{\multirow{2}{*}{NP-hard}} & \multirow{2}{*}{$2$}& \multirow{2}{*}{$\mathbf{2}$}\\
&\multicolumn{2}{c|}{}&\\
&\multicolumn{2}{c|}{\em \cite{Hard-Bottleneck}}&{\em Prop.~\ref{th: approx-partial_pcenter}}\\
   \hline
\end{tabular}
\caption{Main hardness and approximation results in the paper appear in bold characters.}
\label{table-result}
\end{table}



\section{The {Probabilistic \em{p}-Center} problem {under Pressure}}\label{sec:pb}

\subsection{Definition of the problem}

Let $G$ be an edge-weighted graph; we will denote it $G=(V,E,L)$ with $L =\left(\ell_{ij}\right)_{i,j\in V}$ the matrix of lengths. {If $\mathbb{Q}$ denotes the set of rational numbers,  $L$ has entries in $\mathbb{Q}\cup\{\infty\}$ such that $\ell_{ij}<\infty \Leftrightarrow (i,j)\in E$}. We will denote $\ell_m$ and $\ell_M$, respectively the smallest and the largest edge lengths {(i.e., $\ell_M$ is the largest finite entry in $L$)}.
We will refer as the {\em uniform case} the case where all edge lengths are equal. For all problems we consider in this paper, the objective value is linear with respect to the lengths and feasibility conditions due not depend on the lengths. As a consequence, the uniform case is equivalent to the case where all edge lengths are equal to~1. When dealing with the uniform case we will omit $L$ in the instance.
A {\em mixed graph} is a graph with both directed and non-directed edges. When no ambiguity occurs we will use similar notations for graphs and mixed graphs. In the mixed case, we will just identify directed edges an denote them with an arrow in the related  drawing. All non-directed notions in graphs also apply to mixed graphs by considering the {\em non-directed version} of the mixed graph obtained by replacing directed edges by non-directed ones. Similarly, all directed notions apply to mixed graphs since a mixed graph can be seen as a digraph with non-directed edges replaced by two directed edges in opposite directions.
For instance, when speaking about distances in a mixed graph, paths are meant to respect the edge orientations and thus, the matrix of distances is not symmetric anymore.
In an edge-weighted  graph $G=(V,E,L)$ and two vertices $i,j$, $d(i,j)=+\infty$ if $i$ and $j$ are in different connected components. In a mixed graph, we may have $d(i,j)=+\infty$ with $i$ and $j$ in the same connected component. It just means that there is no path from $i$ to $j$ respecting the orientation of directed edges. For example, in the mixed graph represented in Figure~\ref{fig:evacuation_strategy}, $d(2,6)= 5$ while $d(6,2)=\infty$.

In our motivating application,  fire hazards (or any hazard occurring on vertices) is modeled using scenarios. The landscape is represented by  an edge-weighted  graph $G=(V,E,L)$. 
A scenario is associated with each specific fire outbreak.
We restrict ourselves to single fire outbreak and consequently, each scenario $s$ corresponds to a single vertex $s$ on fire.
This restriction is motivated by our primary focus on a relatively short time period after outbreak which assumes an efficient early warning system.
In this case everybody can escape to a shelter before the fire spreads to adjacent zones.

The {\em operational graph} associated with the scenario $s$, denoted by $G^s$, is a mixed graph $G^s=(V, E^s,L^s)$ obtained from $G=(V,E,L)$ by replacing the edges $(s,v)$ incident to $s$ by directed edges $(s,v)$.  All weights are preserved.
Consequently, in $G^s$, vertex $s$ is no longer accessible from another vertex. 

{For every two vertices $i,j$, the distance from $i$ to $j$ in $G^s$ is denoted $d^s(i,j)$. Note that for all $j \in V \setminus \{ s \}$, we have $d^s(j,s)= + \infty$. }


In this paper, we consider a uniform distribution of probabilities over all scenarios: each scenario $s_i,  i \in V\ $ has probability  $\frac{1}{|V|}$ and these events are all independent.


In most $p$-Center problems under indeterminacy, given a solution $C$ with $p$ vertices or less, and given a scenario $s$, the evacuation distance of a vertex $j$ is usually the shortest distance between $j$ and its nearest shelter, $d(j,C)$. This strategy is not adapted to our context and we consider a different evacuation strategy introduced and explained in our previous paper~\cite{demange2020robust}.
This evacuation strategy induces new evacuation distances to shelters. If $s$ is on fire, we have:
\begin{enumerate}
\item for people on $s$, two cases have to be considered. If a shelter is located on $s$, then people present on vertex $s$ are considered  as safely sheltered in it, otherwise we assume that they first run away from the fire in any direction and after they reach a neighbor $j$, they evacuate to the shortest shelter from $j$ in $G^s$.
\item for people who are not on $s$, say on $j\neq s$, the evacuation distance from $j$ to shelter $k$ corresponds to $d^s(j,k)$ in graph $G^s$, i.e., avoiding vertex $s$.
\end{enumerate}

This evaluation of evacuation distances makes our problem specific compared to the literature and induces some additional complexity.
The justification of this measure for people escaping from $s$ is twofold. First, since the area $s$ may be relatively large, a single scenario may correspond to many possible fire configurations, each prohibiting some paths in the zone. The second motivation is to represent decision under stress, a very important characteristic in emergency management: somebody close to the fire may not take rational decisions when selecting a direction while people in another zone can be assumed to behave more rationally.

For a given set $C\subset V$ seen as a set shelter's locations  and a given scenario $s$, the evacuation distance of a zone $j$ is denoted by $r^s(C,j)$. If a shelter is located on $j$, $r^s(C,j)=0$ otherwise we have:


{\begin{equation}
r^s(C,j)= \left\{ \begin{aligned}
 & d^s(j,C) & & \text{if $j\neq s$} \\
 & \max_{v\in N^+_{G^s}(s)} \{ \ell_{sv} +   d^s(v,C) \} & & \text{otherwise}
\end{aligned}\right.
\label{equ:r^s_j(C)}
\end{equation}}
where $N^+_{G^s}(s)$ is the set of all vertices $v$ such that $(s,v) \in E^s$.

Notice that $r^s(C,j)$ is equal to $+ \infty$ if $j$ can't reach any shelter in $G^s$. 

The  {\em evacuation radius} associated with scenario $s$ is defined as $r^s(C)=\max_{j \in V} r^s(C,j)$. Note that $r^s(C)$ is not equal to the usual radius computed in $G^s$: $r^s(C)\geq \max\limits_{v\in V}d^s(v,C)$.

\begin{example}
This example is adapted from~\cite{demange2020robust} and allows to better understand the {evacuation radius} $r^s(C)$ and the operational graph.
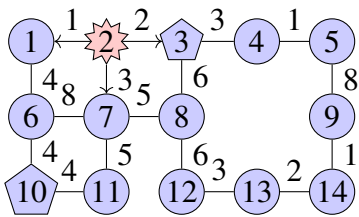
\begin{figure}[hbt]
\centering
 \begin{tikzpicture}

	\node[sommet] (v1) at (0,0){$1$} ;
	\node[sommetFeu] (v2) at (1,0){$2$} ;
	\node[centre] (v3) at (2,0){$3$} ;
	\node[sommet] (v4) at (3,0){$4$} ;
	\node[sommet] (v5) at (4,0){$5$} ;
	
	\node[sommet] (v6) at (0,-1){$6$} ;
	\node[sommet] (v7) at (1,-1){$7$} ;
	\node[sommet] (v8) at (2,-1){$8$} ;
	\node[sommet] (v9) at (4,-1){$9$} ;
	
	\node[centre] (v10) at (0,-2){$10$} ;
	\node[sommet] (v11) at (1,-2){$11$} ;
	\node[sommet] (v12) at (2,-2){$12$} ;
	\node[sommet] (v13) at (3,-2){$13$} ;
	\node[sommet] (v14) at (4,-2){$14$} ;

  	
  	\draw [->] (v2) edge node[above]{1} (v1)  ;
  	\draw [->] (v2) edge node[above]{2} (v3)  ;
  	\draw [->] (v2) edge node[right]{3} (v7)  ;
  	
  	\draw (v1) edge node[right]{4} (v6)  ;
  	\draw (v3) edge node[above]{3} (v4)  ;
  	\draw (v3) edge node[right]{6} (v8)  ;
  	\draw (v4) edge node[above]{1} (v5)  ;
  	\draw (v5) edge node[right]{8} (v9)  ;
  	
  	\draw (v6) edge node[above]{8} (v7)  ;
  	\draw (v6) edge node[right]{4} (v10)  ;
  	\draw (v7) edge node[above]{5} (v8)  ;
  	\draw (v7) edge node[right]{5} (v11)  ;
  	\draw (v8) edge node[right]{6} (v12)  ;
  	\draw (v9) edge node[right]{1} (v14)  ;
  	
  	\draw (v10) edge node[above]{4} (v11)  ;
  	\draw (v12) edge node[above]{3} (v13)  ;
  	\draw (v13) edge node[above]{2} (v14)  ;
 \end{tikzpicture}
\caption{The operational graph associated with scenario 2 with  $C=\{3,10\}$}
\label{fig:evacuation_strategy}
\end{figure}

{Let us consider $p=2$ and 
 the non-directed version $G=(V,E,L)$ of the graph  in Figure~\ref{fig:evacuation_strategy}. We consider  the scenario $s=2$. The related operational mixed graph} is given in Figure~\ref{fig:evacuation_strategy}. 
 Vertices $3$ and $10$ represented by pentagons correspond to shelters' locations ($C=\{3,10\}$). In case of fire on vertex $2$ (scenario $2$), the modification of the graph and the evacuation strategy induce:
\begin{itemize}
\item The shortest path length from $1$ to $3$ is no longer $3$ but $23$, using the shortest path $1, 6, 7, 8, 3$. Consequently, the nearest shelter from vertex $1$ is $10$ 
at a distance of $8$. Thus the evacuation distance  of $1$ in scenario $2$, is equal to $8$ and vertex $1$ is evacuated to vertex $10$.
\item To compute the evacuation distance of vertex $2$ in scenario $2$, we have to consider three neighbors:
\begin{itemize}
\item for neighbor $1$, the distance to the nearest shelter $10$ is 1 + 8 = 9;
\item for neighbor $7$, the distance to the nearest shelter $10$ is 3 + 9 = 12;
\item for neighbor $3$ with a shelter, the distance is 2.
\end{itemize}
Consequently, $r^{2}(C,2)=12$.
\item The {evacuation radius} of the scenario 2 is given by $r^{2}(C)=\max_{j=1, \ldots, 14}r^{2}(C,j)= r^{2}(C,13) = 15$.
\end{itemize}
\qed
\end{example}

{We are now ready to define the problem {\tt Min P$p$CP}. A {\tt Min P$p$CP}-instance will be an edge-weighted graph $G$ and an integer $p$  and  a  solution $C$ will correspond to a set of at most $p$ vertices where to locate  shelters}.
A given solution $C$ corresponds to $n=|V|$ {evacuation radius} $r^1(C), \ldots, r^n(C)$ for $n$ different scenarios. We associate to $C$  the expected value ${\mathbb E}(C)$ of these {evacuation radius} over all scenarios:

\begin{equation}\label{Eq:defE}
 {\mathbb E}(C)= \frac{1}{|V|} \sum\limits_{s \in V}  r^s(C) = \frac{1}{|V|} \sum\limits_{s \in V}  \max_{j\in V} r^s(C,j)
\end{equation}

{$E(C)$ is called probabilistic radius}.
For any set $C$ of centers, {it} can be computed in polynomial time: for each scenario it requires to compute the matrix of shortest path values in the related operational graph, which requires $O(|V|^3)$ operations. So, ${\mathbb E}(C)$ can be computed in $O(|V|^3)$.
The {\tt Min P$p$CP} problem is then to determine a solution $C^*$ minimizing $E$.

We synthesize below the formal definition of the problem:

\begin{problem} 
 \problemtitle{{\tt Min P$p$CP}}
 \probleminput{An edge-weighted graph $G=(V,E,L)$ and an integer $p$ ; the instance is denoted $(G,p)$}
 \problemfeasible{Any $p$-center $C\subset V, |C|\leq p$ satisfying ${\mathbb E}(C)<\infty$ (see Relation~\ref{Eq:defE})}
 \problemoptimization{Minimizes ${\mathbb E}(C)$. }
\end{problem}

{In a more general setting we could add a probability distribution on vertices but in this work we only consider the uniform probability distribution. In this context, recall that by {uniform}, we mean $L$ is the matrix $(\ell_{ij})$ with $\ell_{ij}=1 \Leftrightarrow (i,j)\in E$ and $\ell_{ij}=\infty$ otherwise.}

Note that in our definition, $p$ is part of the instance. We can define natural sub-problems by restricting the possible values for $p$. If $p$ is a fixed value, then the related sub-problem is polynomial since all possible $p$-centers can be enumerated in polynomial time and the 
{probabilistic radius} (objective value) of each one can be determined in polynomial time.

\subsection{Related work}

\cec{
A variant of the {\tt Min $p$-Center} problem for large-scale emergencies is proposed in~\cite{HKM10}, where the disaster affects a single vertex $s$, including any facility on this vertex. This model incorporates both indeterminacy in the facility availability and in the demand: any facility on an affected vertex is no longer available and only the population on this vertex requires evacuation. Our context is really different since we consider that all zones must be evacuated in each scenario $s$ and that a shelter always secures at least the people from the corresponding area.}

\cec{Numerous models for {\tt Min $p$-Center} \mar{under indeterminacy} have already been developed. 
In a non-deterministic environment,  problems are generally described in two stages: first, before the indeterminacy is resolved \mar{(i.e., when the instance is still subject to indeterminacy)}, we \mar{need to} choose locations. \mar{Then, once 
the effective instance is known}, we can react, for example by assigning vertices to 
\mar{centers}.
This description matches real life situations like facility breakdown or natural disaster cutting off communication, and it has been \mar{addressed using various approaches} 
(see for example the reviews~\cite{CNP12, LNS15, S06, CS15}). We \mar{briefly} present some {\tt Min $p$-Center} variants under indeterminacy, as well as some other models relevant for our context.}

\cec{In some models, non determinate parameters may vary independently one from each other, for example in~\cite{AB97, A03, L13, TS12} 
\mar{lengths in the \had{graph} are} described 
\mar{as intervals}. 
In our context, this independence hypothesis is not relevant since, if a fire ignites on a vertex, then all lengths of the edges incident to this vertex are modified in a same way. So, we 
\mar{focus on} the {\tt Min $p$-Center} variants where indeterminacy is represented by a set of discrete scenarios.}

\cec{
\mar{In such a} decision-making environment, we \mar{usually} distinguish two contexts: \mar{risk and uncertainty.} 
\mar{Risk refers to situations where the} values of some parameters are governed by given probability distributions. 
In uncertainty \mar{on the contrary}, 
no probabilistic information is used, either because it is not available or because the decision maker prefers not to resort to it. The first context is referred as {\em stochastic}, or {\em probabilistic} models and the second one corresponds to {\em robust} models. In this second context, a 
measure \mar{of robustness} is usually considered for evaluating the performance of a solution. For example, in~\cite{DZL20}, 
a robust variant of {\tt Min $p$-Center} with scenarios \mar{is studied}. Here, we assume 
\mar{some probability distribution over} scenarios; for this reason, 
we focus only on stochastic/probabilistic variants of {\tt Min $p$-Center} with a \mar{fixed} set of discrete scenarios.}

\cec{In stochastic optimization, one generally 
optimizes the expected value of a given objective function or one maximizes the probability that the solution is ``good''. 
\mar{Such} problems can be solved  using a 
\mar{specific} algorithm, like in~\cite{MAR17}, or  using the general stochastic programming techniques, as in~\cite{BY18}. The problem studied in~\cite{MAR17} is different from ours since the indeterminacy is associated to the demand and it is not possible to re-affect the vertices to a 
\mar{center}. Note that, \mar{despite} {\tt Min P$p$CP} falls into the paradigm of stochastic optimization, 
we cannot 
\mar{easily reduce it to classical $p$-center variants. In particular}, {\tt Min P$p$CP} cannot be seen as a variant in which new distances are associated to each scenario. Indeed, due to our specific evacuation \mar{process, in particular for people in the zone on fire,  the new} 
\had{evacuation} distances \mar{to a shelter} do not systematically correspond to a shortest path in \mar{the new graph} $G^s$. }

Finally, {\tt Min P$p$CP} has been introduced in ~\cite{DHM18}. 
\mar{We have proposed} an explicit solution for the uniform case (all edge lengths are~1) on paths and cycles. In these cases, a solution is characterized by the list of lengths of  segments between two consecutive centers. A $p$-center is then called {\em balanced} if the maximum difference between two segment lengths is minimized and it is {\em monotone} if the sequence of segment lengths is monotone. It is straightforward to show that a balanced solution is optimal for the usual {\tt Min $p$-Center} and, in~\cite{DHM18} we have shown that a monotone balanced solution is also optimal for {\tt Min P$p$CP}. Even though the result is not surprising, the proof was surprisingly non-trivial. In~\cite{DHM18}, we proposed as well some related hardness results. In particular, we showed that {\tt Min P$p$CP} is not approximable on planar graphs of degree~2 or~3 within a ratio less than $\frac{20}{19}$, unless P=NP.  Refining this result in restricted classes of subgrids and designing approximation algorithms were 
\mar{outlined as} open questions.

Note finally that in~\cite{demange2020robust}, we investigated a variant, called {\em robust}, where the objective is to minimize the maximum (worst) evacuation radius over all scenarios instead of minimizing their expected value. It falls into the uncertainty paradigm. For this version, we proposed NP-hardness results in various classes of graphs that include subgrids. Our application motivates this class. We also proposed exact algorithms based on Integer Linear Programming formulation.

In the next subsection we characterize the set of feasible solutions of {\tt Min P$p$CP}.

\subsection{Feasible solutions}
\label{ssec:feasible_solutions}

In this subsection we analyze necessary and sufficient conditions for a solution to be feasible  for a given {\tt Min P$p$CP}-instance $(G,p)$.
Without loss of generality we will consider that $G=(V,E,L)$ is a connected graph.
A vertex $a\in G$ is an {\em articulation point} if and only if removing $a$ disconnects the graph $G$.
We denote by $\mathcal{A}(G)$ the set of articulation points of $G$.

We call {\em articulation component}  of $G$ associated with an articulation point $a$  a connected component of $G\setminus\{a\}$.
Then every vertex $a\in \mathcal{A}(G)$ is associated to at least 2 articulation components, and every articulation component is associated to one articulation point. A graph is 2-connected if it has no articulation point; in this case there is no articulation component.

A  {\em minimal articulation component}, or MAC for short, is an articulation component that does not strictly contain another articulation component.
We denote $\Upsilon(G)$ the set of minimal articulation components.
Note that an articulation component that is a singleton $\{v\}$ is necessarily minimal and this occurs if and only if $v$ is  a vertex of degree~1.

\begin{lem}
\label{lem:feasibility_1}
$A$ is a minimal articulation component of $G$ if and only if $A$ is an articulation component which does not include an articulation point of $G$.
\end{lem}
\begin{pf}

$\Rightarrow$ By contrapositive we prove that if an articulation component $A$ includes an articulation point, then $A$ is not minimal.
Let $A$ be an articulation component induced by the articulation point $a\in V$. Suppose $b\in A$ is an articulation point of $G$.
Then $b$ induces at least two disjoint connected components in $G \setminus \{b\}$.
Since  $b\neq a$, $a$ is in  one connected component of $G \setminus \{b\}$, consequently $G\setminus A$ is a subset of this connected component.
It follows that at least another component of $G \setminus \{b\}$ is {contained} in $A$, which means that $A$ is not minimal.

$\Leftarrow$ The proof is also by contrapositive.
We prove that if $A$ is a non-minimal articulation component, then $A$ includes an articulation point.
Let $A$ an articulation component that is not minimal.
Then there is an articulation component $B\subset A$ induced by the articulation point $b\in V$, such that $B\neq A$.
Consider $x\in A\setminus B$ and $y\in B$. Since $A$ is connected,  $x$ and $y$ are connected in $A$ by a path; this path necessarily  crosses $b$ and in particular  $b\in A$.
\qed
\end{pf}

\begin{lem}
\label{lem:feasibility_2}
All minimal articulation components of $G$ are pairwise disjoints.
\end{lem}
\begin{pf}
By contrapositive, we assume $A\in \Upsilon(G)$ and $B$ an articulation component such that $B\neq A$ and $B\cap A \neq \emptyset$.
We prove then that $B$ is not minimal.

Let $x\in A\cap B$.
Since $A$ is a MAC, $B \not\subset A$.
Then there is a vertex $y\in B\cap (V\setminus A)$.
Every path between $x$ and $y$ in $G$ crosses $a$.
As $B$ is a connected component, there is a path from $x$ to $y$ in $B$, thus $a\in B$, and $B$ is not a MAC by Lemma~\ref{lem:feasibility_1}.  \qed
\end{pf}

Given an edge-weighted graph $G=(V,E,L)$ and $p$, we denote with $\mathcal{C}_p(G)$ the set of feasible solutions of the {\tt Min P$p$CP}-instance $(G,p)$.

\begin{proposition}
\label{prop:feasibility_condition}
Let $(G,p)$ be an instance of {\tt Min P$p$CP} with $|V|\ge 2$. A solution $C\subset V, |C|\leq p$  is in $\mathcal{C}_p(G)$ if and only if $|C|\ge 2$ and $C$ includes at least one vertex in each minimal articulation component of $G$.
\end{proposition}
\begin{pf}
Suppose $C$ is a feasible solution for {\tt Min P$p$CP} on $G$.
We have seen that $C$ is a feasible solution for {\tt Min P$p$CP} if and only if $r^s(C,j) \in \mathbb{R}, \forall j,s\in V$, i.e. all the evacuation distances over all vertices and all scenarios are finite.

First suppose there is no articulation point, then $G$ has no articulation components.
Let $s\in C$, and $x\in V, x\neq s$. In scenario $s$, $x$ is assigned to a center that is not $s$.
Thus $|C|\ge 2$.
Conversely, if $|C|\ge 2$, for any scenario $s$, $G\setminus {\{s\}}$ is connected and contains at least one center.

Second, suppose $G$ has at least one articulation point and consequently at least 2 disjoint articulation components.
In addition, if $A$ is an articulation component of $G$ induced by the articulation point $a$, then $\forall j\in A, r^a(C,j) \in \mathbb{R}$ if and only if $C\cap A \neq \emptyset$.
Then $C$ intersects all articulation components. In particular $|C|\ge 2$ and $C$ intersects all minimal ones.
Conversely, if $C$ intersects all MACs then $|C|\ge 2$ and it intersects all articulation components since any articulation component contains a MAC.
\qed
\end{pf}

\begin{remarque}
{Feasibility  is weight-independent.}
\end{remarque}

\begin{corollary}If $G$ has at least 2 vertices, $\mathcal{C}_1(G)= \emptyset$.
\end{corollary}
As a consequence, from now we will consider only {\tt Min P$p$CP} instances satisfying $p\geq 2$.

\begin{corollary}For a given $p$, we can verify in polynomial time whether $\mathcal{C}_p(G)\neq \emptyset$.
\end{corollary}
\begin{pf}
For $G=(V,E)$, we generate $\mathcal{A}(G)$ in $O(|V|+|E|)$ using Tarjan's Algorithm (\cite{tarjan1972depth}).
The minimal connected components of $G$ are the connected components of $G\setminus\mathcal{A}(G)$ adjacent to at most one articulation point in $G$, where a set $V'$ of vertices is said adjacent to a vertex if this vertex has at least one neighbor in $V'$.

There is a feasible solution for {\tt Min P$p$CP} on $G$ if $p$ is greater or equal to the number of MACs.
\qed \end{pf}


\begin{corollary}
For all $C\in\mathcal{C}_p(G)$, $C$ necessarily includes all vertices of degree~1.
\end{corollary}
\begin{pf}
Every vertex of degree~1 is a MAC of $G$. Then by Proposition~\ref{prop:feasibility_condition}, a feasible solution includes all vertices of degree~1.
\qed\end{pf}

%
%
%
%
%



\subsection{Further notations}\label{subsec:further notations}

For a graph $G=(V,E)$ and a set $V'\subset V$, we will denote $G[V']$ the {\em subgraph} of $G$ induced by $V'$. $G[V']$ is called a subgraph of $G$. A {\em partial graph} of $G$ is a graph $(V,E')$ with $E'\subset E$ obtained from $G$ by deleting 0 or some edges. A partial subgraph of $G$ is a partial graph of a subgraph of $G$. For $U\subset V$, we denote $G\setminus U$ the graph  $G[V\setminus U]$. {A {\em pending} vertex in a graph is a vertex of degree~1. A $n\times m$ {\em grid} is the graph $\mathcal{G}=(\mathcal{V,E})$ with vertex set $\mathcal{V}=\{(i,j), i\in \{0, \ldots n-1\}, j\in \{0, \ldots, m-1\}\}$ and $((i,j),(k,l))\in \mathcal{E}$ if and only if  $|i-k|+|j-l|=1$. A (partial) {\em subgrid} is a (partial) subgraph of a grid. For instance, the graph in Figure~\ref{fig:evacuation_strategy} is a partial subgrid. Given a subgrid $G=(V,E)$, a grid embedding is a one-to-one function from $V$ to  $\mathcal{V}$ for some dimensions $(n,m)$ such that every edge $(u,v)\in E$ maps to an edge of the $n\times m$ grid. If $u\in V$ maps to $(i,j)$ in the grid, $(i,j)$ are called the coordinates of $u$. Unless otherwise stated, each time we will refer to a subgrid, we will assume that a grid embedding is given. As defined in~\cite{color-grid}, for a partial subgrid $G$ and a positive integer $f$, the $f$-expansion of $G$, denoted $Exp(G,f)$, is obtained  from $G$ by inserting $f-1$ vertices on each edge (each edge becomes a path of $f$ edges).  If $f\geq 2$, the {\em $f$-expansion} of any partial subgrid is a subgrid. If $G$ is a subgrid embedded in a $n\times m$  grid $\mathcal{G}$, then $Exp(G,f)$ is a subgrid embedded in the  $[(n-1)f+1]\times [(m-1)f+1]$ grid $Exp(\mathcal{G},f)$. The vertex set of $G$ can be seen as a subset of the vertex set of $Exp(G,f)$ and more precisely, in the related grid embedding of $Exp(G,f)$, the coordinates of any vertex $u\in V$ are multiplied by $f$ compared to its coordinates in the original grid embedding of $G$ in $\mathcal{G}$. Subgrids, and to a lesser extent partial subgrids, constitute a natural class of instances in our motivating application. It corresponds to the case where the landscape is divided into square areas and  some areas are not considered since they correspond for instance to natural barriers, like lakes, or to protected private lands that can neither been used for sheltering nor for evacuating.}

\section{Hardness result}
\label{sec:complexity}

In all this section we assume that all edge lengths are~1 {(uniform case)}.
We remind that  {\tt Min $p$-Center} is NP-hard for $p\ge 2$ on planar graphs of maximum degree at least~3 (\cite{KH79}). This result does not  immediately imply the hardness of {\tt Min P$p$CP}.
Indeed, we defined our model with fixed uniform probabilities, which does not count the classic deterministic {\tt Min $p$-Center} problem as one of its specific cases.

\mar{Note first that the decision version of {\tt Min P$p$CP} is in NP. Indeed, if we consider a set of centers $C$, for each scenario $s$, the evacuation radius $r^s(C)$ can be computed in polynomial time using a shortest path algorithm. Then, ${\mathbb E}(C)$ can be computed in polynomial time using Relation~\ref{Eq:defE}. So, a non-deterministic algorithm will infer a $p$-center and verify in polynomial time whether its value exceeds or not the target. As a consequence, any hardness in approximation result for  {\tt Min P$p$CP} can be immediately turned into a NP-completeness result.}

In~\cite{DHM18}, we showed that {\tt Min P$p$CP} cannot be approximated on planar graphs of degree~2 or~3 with a ratio less than $\frac{20}{19}$. \mar{Actually, a close look on the proof shows  that the hardness result holds for bipartite planar graphs of degree~2 or~3. } In this section, we prove that {\tt Min P$p$CP} cannot be approximated with a ratio less than $\frac{56}{55}$  on a restricted subclass of bipartite planar graphs, the class of subgrids with degree at most~3.
\mar{In terms of NP-completeness, it shows that {\tt Min P$p$CP} is NP-complete in subgrids,  while our previous result  established it for the  class of bipartite planar graphs. }

The proof uses two classical optimization problems in graphs, \mar{the {\tt Min Dominating Set} and the {\tt Min Vertex Cover} problems}. A {\em dominating set} in a graph $G = (V, E)$ is a subset $U$ of $V$ such that every vertex not in $U$ is adjacent to at least one member of $U$. The {\tt Min Dominating Set} problem is to find a dominating set of minimum size. We will denote by $\gamma(G)$ the minimum size of a dominating set in $G$. 
The {\tt Min Dominating Set} problem is shown NP-hard on subgrids in~\cite{CLARK1990165}. Note the relation between  { \tt Min Dominating Set} and the deterministic {\tt Min $p$-Center} problem: for a graph $G=(V,E)$, $U\subseteq V$ is a dominating set, if and only if $U$ is a $|U|$-center of radius~1. \mar{In Subsection~\ref{subsec:relations_PpCP_DS}, we will establish links between { \tt Min Dominating Set} and {\tt Min P$p$CP}.}

A {\em vertex cover} of a graph is a set of vertices such that each edge of the graph is incident to at least one vertex of the set.
The {\tt Min Vertex Cover} problem is to find a vertex cover of minimum size.
We will denote by $\tau(G)$ the minimum size of a vertex cover in $G$.
In~\cite{KH79},  {\tt Min Vertex Cover} is shown NP-hard on planar graphs of maximum degree~3. \mar{It is straightforward  to show that it remains NP-hard on planar graphs with vertices 2 or 3 (see~\cite{DHM18})}.


\mar{The proof requires a technical lemma similar to Lemma~7 in~\cite{DHM18}. Since it is slightly different, we give a proof in Appendix. }

\begin{lem}\label{lem:relation_size_minVC} (\cite{DHM18})\\
Let $G=(V,E)$ be a graph and $G'=(V',E')$ be the graph obtained by inserting $2k_{uv}$ vertices on each edge $(u,v)\in E$, where $k_{uv}$ is a non-negative integer. Then we have
$$ \tau(G') = \tau(G) + \sum_{uv\in E} k_{uv}$$
\end{lem}





The remaining of the section is dedicated to prove Theorem~\ref{thm:complexity}. 
In Subsection~\ref{subsec:blueprint} 
we explain the general scheme of the demonstration before giving all details in Subsections~\ref{subsec:Transformations},~\ref{subsec:relations_PpCP_DS} and~\ref{subsec:theorem_complexity}.


\subsection{Global blueprint of the proof}
\label{subsec:blueprint}



{In Theorem~\ref{thm:complexity}, we will show that a polynomial time approximation algorithm $\mathtt{A}$ for {\tt Min P$p$CP} in subgrid of degrees~$\{2,3\}$ ($p$ being part of the instance)  guaranteeing a ratio of at most $\frac{56}{55}$ could be used  to compute in polynomial time the size of the minimum vertex cover on a planar graph of degrees~$\{2,3\}$, which is a contradiction.}

We adapt the proof of the hardness result in~\cite{DHM18} to obtain a hardness result in subgrids. The requirement that the resulting graph is a subgraph of a grid induces significant technical difficulties for both the reduction and its analysis.

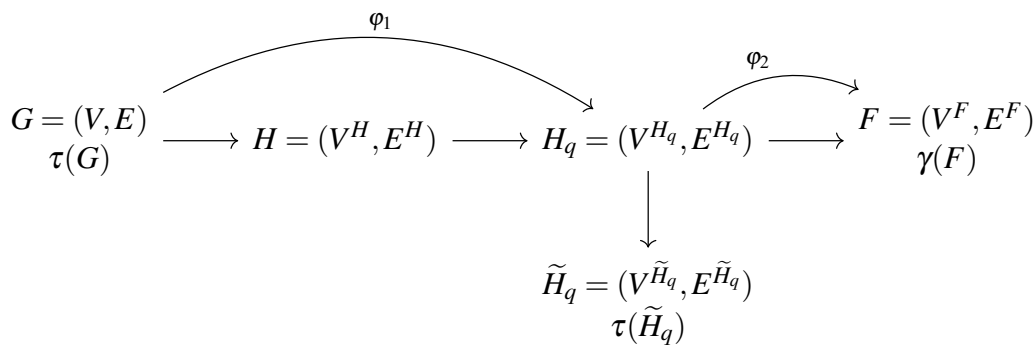
\begin{figure}[hbt]
\centering
\begin{tikzcd}
\begin{tabular}{@{}c@{}}$G=(V,E)$\\$\tau(G)$ \end{tabular} \arrow[r] \arrow[rr,bend left, "\varphi_1"]	& H=(V^H,E^H)\arrow[r] & H_q=(V^{H_q},E^{H_q})\arrow[d]\arrow[r] \arrow[r,bend left, "\varphi_2"]	& \begin{tabular}{@{}c@{}}$F=(V^F,E^F)$\\$\gamma(F)$ \end{tabular} \\
			& &\begin{tabular}{@{}c@{}}$\widetilde H_q =(V^{\widetilde H_q},E^{\widetilde H_q})$\\$\tau(\widetilde H_q)$ \end{tabular} &
\end{tikzcd}
\caption{The different graphs {involved in the reduction}.}
\label{fig:schema_transformation}
\end{figure}

{We start from a  planar graph $G = (V,E)$ with degrees~$\{2,3\}$, instance of {\tt Min Vertex Cover}.
We randomly choose an orientation of the edges of $G$ that will be used in our reductions and analysis. We then apply successively two transformations, Transformation~\ref{transformation:1}, denoted $\varphi_1$ and Transformation~\ref{transformation:3}, denoted $\varphi_2$ that are detailed in Subsection~\ref{subsec:Transformations}.
Figure~\ref{fig:schema_transformation} gives a simple schematic representation of the  whole reduction. }

Transformation~\ref{transformation:1} ($\varphi_1$) constructs from $G$ a subgrid  ${H_q}=(V^{H_q}, E^{H_q})$, for some positive integer $q$ specified later, in such a way that:
\begin{itemize}
    \item[$\bullet$] $V\subset V^{H_q}$,
    \item[$\bullet$] Edges $(u,v)$ of $G$ map   to non-crossing paths $P^{H_q}_{uv}$ of even length between $u$ and $v$ in $H_q$.
\end{itemize}
The subgrid $H$ appearing in Figure~\ref{fig:schema_transformation} is an intermediate stage not directly used in the analysis.

{We then apply Transformation~\ref{transformation:3} ($\varphi_2$) to construct a subgrid $F =(V^F, E^F)$ from $H_q$. Roughly speaking,  it consists in replacing the first two edges of $P^{H_q}_{uv}$ (where $uv\in E$ is oriented from $u$ to $v$) with a gadget $T^2$, and every other edge of $H_q$ with a gadget $T^1$, both defined in the next subsection.}

{For the analysis now, we note that there is no direct and easy link between $\tau(G)$ and $\tau(H_q)$ since $\tau(H_q)$ can be obtained in polynomial time ($H_q$ is bipartite) while $G$ is meant to be an instance of an NP-hard restriction of {\tt Min Vertex Cover}.
For this reason, we introduce an auxiliary graph $\widetilde H_q=(V^{\widetilde H_q},E^{\widetilde H_q})$.
It can be seen as
 a perturbation of $H_q$ with a direct link between $\tau(\widetilde H_q)$ and $\tau(G)$. It is simply obtained  by replacing, for every edge $(u,v)\in E$,  the two first edges of  the path  $P^{H_q}_{uv}$ by a single edge. This way, the path  $P^{H_q}_{uv}$ of even length becomes, in $\widetilde H_q$, a path $P^{\widetilde H_q}_{uv}$ of odd length and  Lemma~\ref{lem:relation_size_minVC} can be used to write $\tau(\widetilde H_q)$ as a function of $\tau(G)$.}

 {On the other hand, as outlined in Lemma~\ref{lem:link_VC_DS}, the properties of the two gadgets allow to establish a direct link between dominating sets in $F$ and vertex covers in $\widetilde H_q$. In all, it gives a relation between the {\tt Min Dominating Set} problem in $F$ and the {\tt Min Vertex Cover} problem in $G$.}



Then, in Subsection~\ref{subsec:relations_PpCP_DS}, we outline different relations between the {\tt Min Dominating Set} problem and {\tt Min P$p$CP} in a triangle-free graphs without pending vertices using three lemmas, Lemma~\ref{lem:relation_DS_Radius}, Lemma~\ref{lem:relation_value_min} and Lemma~\ref{lem:radius}. This can be applied to $F$.  


Finally, in Subsection~\ref{subsec:theorem_complexity}, we use these results to establish  Theorem~\ref{thm:complexity}.
We show that, 
when applying $\mathtt{A}$ on $F$ for $p< \gamma(F)$, the output is a solution of {\tt Min P$p$CP}  of 
{probabilistic radius} at least~2, while applying it for $p=\gamma({F})$ gives a solution of  {probabilistic radius} less than~2.
Hereby we can use such an algorithm to compute  $\gamma(F)$, and consequently $\tau(G)$. Since constructing $H$, $H_q$, $\widetilde H_q$ and $F$, as well as evaluating the value of a {\tt Min P$p$CP} solution, can be done in polynomial time, and since algorithm~$\mathtt{A}$ is applied less than $|V|$ times, the whole process is polynomial.

\subsection{Details on the transformations and their properties}
\label{subsec:Transformations}






\begin{transformation}{From a  planar graph $G=(V,E)$ to a subgrid ${H_q} = (V^{H_q},E^{H_q})$ with $q > 0$.}
\label{transformation:1}
\end{transformation}

Using a result of~\cite{embed}, we can embed $G=(V,E)$ in a grid $H=(V^H,E^H)$ of polynomial size. Vertices of $G$ are mapped to vertices of the grid, and {edges $(u,v)$ of $G$ map to non-crossing paths $P_{uv}^H$ between $u$ and $v$ in the grid. Note that we cannot control the length and parity of these paths. The resulting graph is a partial subgrid and not necessarily a subgrid yet. We then perform a $2q$-expansion for some positive integer $q$ specified later. The resulting graph   ${H_q}=(V^{H_q}, E^{H_q})$ is a subgrid ($q>0$). In addition, since the expansion multiplies by $2q$ all path lengths from $H$ to $H_q$, edges $(u,v)$ of $G$ map to non-crossing paths $P_{uv}^{H_q}$ of even length between $u$ and $v$ in $H_q$. It means that paths $P_{uv}^{H_q}$ have $2k_{uv}+1$ internal vertices (excluding $u$ and $v$) for some non-negative integers $k_{uv}$. }

\begin{figure}[hbt]
\centering
\begin{tikzpicture}

\node[sommet] (a) at (5,0){a} ;
\node[sommet] (b) at (5,4){b} ;
\node[sommet] (c) at (5,8){c} ;
\node[sommet] (d) at (1,4){d} ;

\node[sommet] (ad) at (1,0){} ;
\node[petitsommet] () at (2,0){} ;
\node[petitsommet] () at (3,0){} ;
\node[petitsommet] () at (7,0){} ;
\node[petitsommet] () at (8,0){} ;
\node[sommet] (ac) at (9,0){} ;

\node[petitsommet] () at (1,1){} ;
\node[petitsommet] () at (9,1){} ;
\node[petitsommet] () at (1,2){} ;
\node[petitsommet] () at (5,2){} ;
\node[petitsommet] () at (9,2){} ;
\node[petitsommet] () at (1,3){} ;
\node[petitsommet] () at (5,3){} ;
\node[petitsommet] () at (9,3){} ;

\node[petitsommet] () at (1,5){} ;
\node[petitsommet] () at (9,5){} ;
\node[petitsommet] () at (1,6){} ;
\node[petitsommet] () at (5,6){} ;
\node[petitsommet] () at (9,6){} ;
\node[petitsommet] () at (1,7){} ;
\node[petitsommet] () at (5,7){} ;
\node[petitsommet] () at (9,7){} ;

\node[petitsommet] () at (2,4){} ;
\node[petitsommet] () at (3,4){} ;
\node[sommet] () at (9,4){} ;

\node[sommet] (cd) at (1,8){} ;
\node[petitsommet] () at (2,8){} ;
\node[petitsommet] () at (3,8){} ;
\node[petitsommet] () at (6,8){} ;
\node[petitsommet] () at (7,8){} ;
\node[petitsommet] () at (8,8){} ;
\node[sommet] (ca) at (9,8){} ;

\begin{scope}[on background layer]
\foreach \from/\to in {a/b,b/c,b/d,a/ad,a/ac,ad/d,ac/ca,d/cd,cd/c,c/ca}
   	\draw (\from) -- (\to) ;
\end{scope}

\node[align=center,text width=4cm] at (5,-1) {$\widetilde H_q$};

\end{tikzpicture}
\caption{The graph $\widetilde H_q$ obtained from $G$ through ${H_q}$.\\
$k_{ab}=k_{bc}=k_{bd}=1$, $k_{ad}=k_{cb}=3$ and $k_{ac}=7$}
\label{fig:transformation2}
\end{figure}
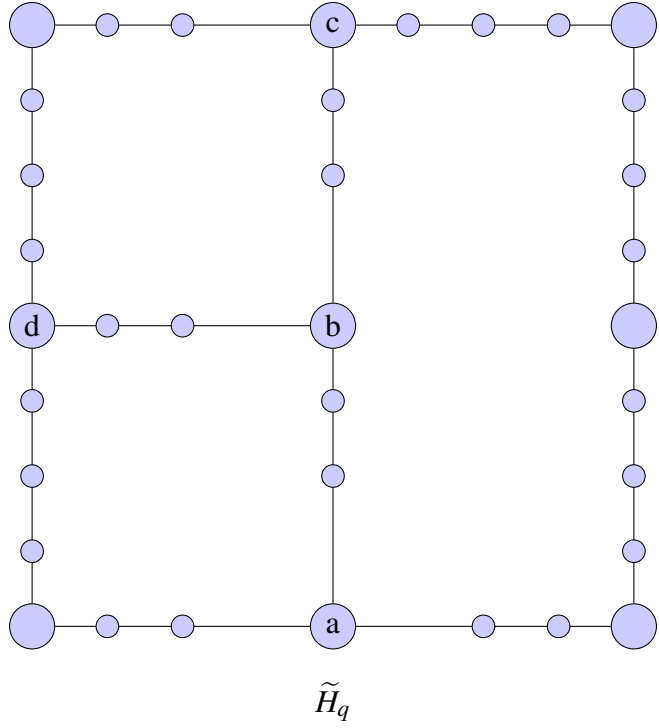

\begin{example}
\label{example:Transformation1}
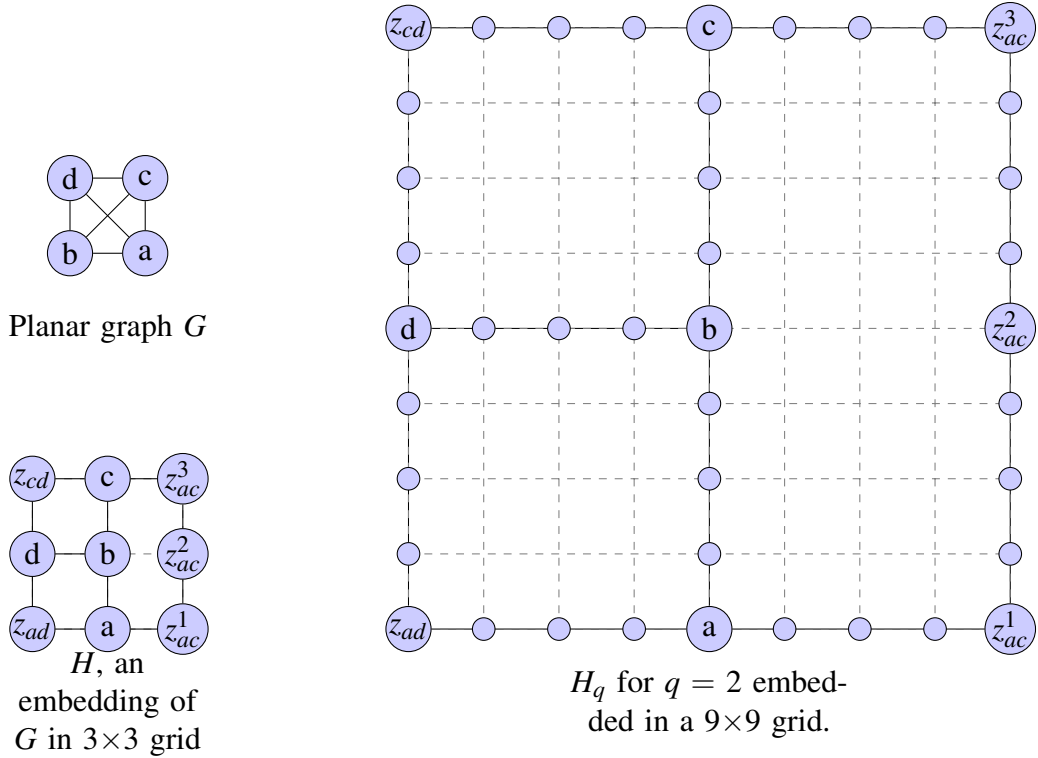
\begin{figure}[hbt]
\begin{tikzpicture}

\node[sommet] (a) at (-2.5,5){a} ;
\node[sommet] (b) at (-3.5,5){b} ;
\node[sommet] (c) at (-2.5,6){c} ;
\node[sommet] (d) at (-3.5,6){d} ;

\foreach \from/\to in {a/b,a/d,b/d,b/c,c/d,a/c}
   	\draw (\from) -- (\to) ;
   	
\node[align=center,text width=4cm] at (-3,4) {Planar graph $G$};

\draw [step=1, dashed, gray] (-4,0) grid (-2,2);


\node[sommet] (a2) at (-3,0){a} ;
\node[sommet] (b2) at (-3,1){b} ;
\node[sommet] (c2) at (-3,2){c} ;
\node[sommet] (d2) at (-4,1){d} ;

\node[sommet] (z1) at (-4,0){$z_{ad}$} ;
\node[sommet] (z2) at (-4,2){$z_{cd}$} ;
\node[sommet] (z3) at (-2,0){$z_{ac}^1$} ;
\node[sommet] (z4) at (-2,1){$z_{ac}^2$} ;
\node[sommet] (z5) at (-2,2){$z_{ac}^3$} ;

\foreach \from/\to in {a2/b2,b2/d2,b2/c2,a2/z1,z1/d2,z2/d2,z2/c2,z5/c2,z3/z4,z4/z5,z3/a2}
   	\draw (\from) -- (\to) ;
   	
\node[align=center,text width=2.5cm] at (-3,-1) {$H$, an embedding of $G$ in 3$\times$3 grid};

\draw [step=1, dashed, gray] (1,0) grid (9,8);

\node[sommet] (a3) at (5,0){a} ;
\node[sommet] (b3) at (5,4){b} ;
\node[sommet] (c3) at (5,8){c} ;
\node[sommet] (d3) at (1,4){d} ;

\node[sommet] (ad) at (1,0){$z_{ad}$} ;
\node[petitsommet] () at (2,0){} ;
\node[petitsommet] () at (3,0){} ;
\node[petitsommet] () at (4,0){} ;
\node[petitsommet] () at (6,0){} ;
\node[petitsommet] () at (7,0){} ;
\node[petitsommet] () at (8,0){} ;
\node[sommet] (ac) at (9,0){$z_{ac}^1$} ;

\node[petitsommet] () at (1,1){} ;
\node[petitsommet] () at (5,1){} ;
\node[petitsommet] () at (9,1){} ;
\node[petitsommet] () at (1,2){} ;
\node[petitsommet] () at (5,2){} ;
\node[petitsommet] () at (9,2){} ;
\node[petitsommet] () at (1,3){} ;
\node[petitsommet] () at (5,3){} ;
\node[petitsommet] () at (9,3){} ;

\node[petitsommet] () at (1,5){} ;
\node[petitsommet] () at (5,5){} ;
\node[petitsommet] () at (9,5){} ;
\node[petitsommet] () at (1,6){} ;
\node[petitsommet] () at (5,6){} ;
\node[petitsommet] () at (9,6){} ;
\node[petitsommet] () at (1,7){} ;
\node[petitsommet] () at (5,7){} ;
\node[petitsommet] () at (9,7){} ;

\node[petitsommet] () at (2,4){} ;
\node[petitsommet] () at (3,4){} ;
\node[petitsommet] () at (4,4){} ;
\node[sommet] () at (9,4){$z_{ac}^2$} ;

\node[sommet] (cd) at (1,8){$z_{cd}$} ;
\node[petitsommet] () at (2,8){} ;
\node[petitsommet] () at (3,8){} ;
\node[petitsommet] () at (4,8){} ;
\node[petitsommet] () at (6,8){} ;
\node[petitsommet] () at (7,8){} ;
\node[petitsommet] () at (8,8){} ;
\node[sommet] (ca) at (9,8){$z_{ac}^3$} ;

\begin{scope}[on background layer]
\foreach \from/\to in {a3/b3,a3/ad,a3/ac,ad/cd,ac/ca,b3/d3,b3/c3,cd/ca}
   	\draw (\from) -- (\to) ;
\end{scope}

\node[align=center,text width=6cm] at (5,-1) {${H_q}$ for $q=2$ embedded in a 9$\times$9 grid.};

\end{tikzpicture}
\caption{Example of Transformation~1}
\label{fig:transformation1}
\end{figure}

Suppose the  planar graph $G=(V,E)$ is a complete graph on four vertices $\{a,b,c,d\}$ as presented in Figure~\ref{fig:transformation1} and set $q=2$.
We choose an orientation of $G$ 
such that the oriented edges of $G$ are $\{(a,b),(a,c),(b,c),(c,d)\}$.
$H=(V^H,E^H)$ corresponds to a possible embedding of $G$ in a grid, where the edge $(a,d) \in E$ maps to the path $\{a,z_{ad},d\}$ in $H$.
Next, we construct the subgrid ${H_q}$ by applying the $2q$-expansion.  The resulting graph ${H_q}$ can be seen on the right side of Figure~\ref{fig:transformation1}. Finally, the related graph $\widetilde H_q$ is represented in Figure~\ref{fig:transformation2}.

\end{example}\qed

{As already noticed in Subsection~\ref{subsec:blueprint}, we cannot establish a direct link between $\tau(G)$ and $\tau(H_q)$ but since we now control the parity of paths $P_{uv}^{H_q}$, it is easy to slightly modify $H_q$ so as we can apply Lemma~\ref{lem:relation_size_minVC}. This is the role of the graph $\widetilde H_q=(V^{\widetilde H_q},E^{\widetilde H_q})$. Recall that this graph is obtained from
$H_q$  by replacing, for every edge $(u,v)\in E$,  the two first edges of  the path  $P^{H_q}_{uv}$ by a single edge, as illustrated in Figure~\ref{fig:transformation2}. This way,
$\widetilde H_q$ can directly be obtained from $G$ by inserting $2k_{uv}$ vertices on each edge $(u,v)\in E$. As a consequence,
Lemma~\ref{lem:relation_size_minVC} allows to establish:}

\begin{equation}\label{eq:tau_G_wideH}
    \tau(\widetilde H_q) = \tau(G) + \sum_{(u,v)\in E} k_{uv}.
\end{equation}

In addition, we have:

\begin{equation}\label{eq:metricsHtilde}\begin{array}{lrl}
|V^{\widetilde H_q}|&=&|V|+2\sum\limits_{e\in E} k_e\\
|E^{\widetilde H_q}|&=&|E|+2\sum\limits_{e\in E} k_e
\end{array}\end{equation}

By construction, we have $\forall (u,v) \in E, 2k_{uv}+1\ge 2q-1$, which gives:

\begin{equation}\label{eq:q}
    \forall (u,v) \in E, k_{uv}\ge q-1.
\end{equation}

\label{example:Transformation2}

\begin{transformation}{From subgrid ${H_q}=(V^{H_q},E^{H_q})$ to subgrid $F=( V^F, E^F)$.  }
\label{transformation:3}
\end{transformation}



{Thanks to the $2q$-expansion,  for $(u,v)\in E$ oriented from $u$ to $v$, the first two edges of $P^{H_q}_{uv}$ in $H_q$ are  both horizontal or vertical. Note as well that the orientation of $G$ immediately defines an orientation of $H_q$ and of $\widetilde H_q$.
}
We can then construct the subgrid $F=(V^F, E^F)$ from the subgrid $H_q$ as follows.

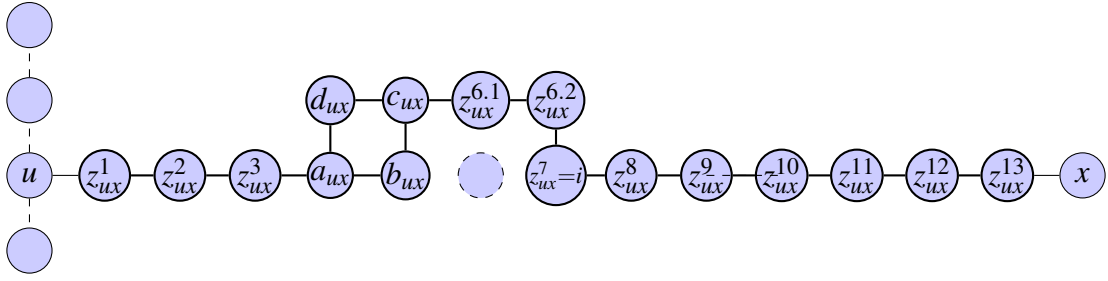
\begin{figure}[hbt]
\centering
  \begin{tikzpicture}
	
	\node[sommet] (n0) at (0,0){${u}$} ;
	\node[sommet] (n17) at (14,0){${x}$} ;
	
	\node[sommet, thick](n1) at (1,0){$z_{ux}^1$} ;
    \node[sommet, thick] (n2) at (2,0){$z_{ux}^2$} ;
	\node[sommet, thick] (n3) at (3,0){$z_{ux}^3$} ;
	\node[sommet, thick] (n4) at (4,0){$a_{ux}$} ;	
	\node[sommet, thick] (n5) at (5,0){$b_{ux}$} ;
    \node[sommet, thick] (n6) at (5,1){$c_{ux}$} ;	
    \node[sommet, thick] (n7) at (4,1){$d_{ux}$} ;
    \node[sommet, thick] (n8) at (6,1){$z_{ux}^{6.1}$} ;
    \node[sommet, thick] (n9) at (7,1){$z_{ux}^{6.2}$} ;
    \node[sommet, thick] (n10) at (7,0){$\scriptstyle{z_{ux}^7=i}$} ;
    \node[sommet, thick] (n11) at (8,0){$z_{ux}^8$} ;

    \node[sommet, thick] (n12) at (9,0){$z_{ux}^9$} ;
    \node[sommet, thick] (n13) at (10,0){$z_{ux}^{10}$} ;
    \node[sommet, thick] (n14) at (11,0){$z_{ux}^{11}$} ;
    \node[sommet, thick] (n15) at (12,0){$z_{ux}^{12}$} ;

    \node[sommet, thick] (n16) at (13,0){$z_{ux}^{13}$} ;

    \node[sommet]  (n18) at (0,1){} ;
    \node[sommet]  (n19) at (0,2){} ;
    \node[sommet]  (n20) at (0,-1){} ;

    \node[sommet,dashed]   at (6,0){} ;


\foreach \from/\to in  {n1/n2,n2/n3,n3/n4,n4/n5,n4/n7,n6/n7,n5/n6, n6/n8,n8/n9,n9/n10, n10/n11,n11/n12,n12/n13,n13/n14,n14/n15,n15/n16}
     	\draw[thick] (\from) -- (\to)  ;
   \foreach \from/\to in  {n0/n1,n16/n17}
    	\draw[] (\from) -- (\to) ;
\foreach \from/\to in  {n0/n18,n18/n19,n0/n20}
    	\draw[dashed] (\from) -- (\to) ;
\draw[dashed] (9,0) -- (10,0);
  \end{tikzpicture}
  \caption{Gadget $T^2_{ux}$ used in $F$ for $(u,i),(i,x)\in E^{ H_q}$; $z_{ux}^7=i$.}
  \label{fig:gadget2}
\end{figure}

\begin{figure}[hbt]
\centering
  \begin{tikzpicture}
	
	\node[sommet][] (n0) at (0,0){${x}$} ;
	\node[sommet](n14) at (7,0){${y}$} ;
	\node[sommet, thick] (n1) at (1,0){$z_{xy}^1$} ;
    \node[sommet, thick] (n2) at (2,0){$z_{xy}^2$} ;
	\node[sommet, thick] (n3) at (2,1){$z_{xy}^3$} ;
	\node[sommet, thick] (n4) at (2,2){$a_{xy}$} ;	
	\node[sommet, thick] (n5) at (2,3){$b_{xy}$} ;
    \node[sommet, thick] (n6) at (3,3){$c_{xy}$} ;	
    \node[sommet, thick] (n7) at (3,2){$d_{xy}$} ;
    \node[sommet, thick] (n8) at (4,3){$z_{xy}^{4.1}$} ;
    \node[sommet, thick] (n9) at (5,3){$z_{xy}^{4.2}$} ;
    \node[sommet, thick] (n10) at (5,2){$z_{xy}^{4.3}$} ;
    \node[sommet, thick]  (n11) at (5,1){$z_{xy}^{4.4}$} ;
    \node[sommet, thick]  (n12) at (5,0){$z_{xy}^{5}$} ;
    \node[sommet, thick]  (n13) at (6,0){$z_{xy}^{6}$} ;


    \node[sommet,dashed]   at (3,0){} ;
    \node[sommet,dashed]  at (4,0){} ;

\foreach \from/\to in  {n1/n2,n2/n3,n3/n4,n4/n5,n4/n7,n6/n7,n5/n6, n6/n8,n8/n9,n9/n10,n10/n11, n11/n12,n12/n13}
     	\draw[thick] (\from) -- (\to)  ;
   \foreach \from/\to in  {n0/n1,n13/n14}
    	\draw[] (\from) -- (\to) ;
  \end{tikzpicture}
  \caption{Gadget $T^1_{xy}$ used in $F$ for $(x,y)\in E^{H_q}$.}
  \label{fig:gadget1}
\end{figure}
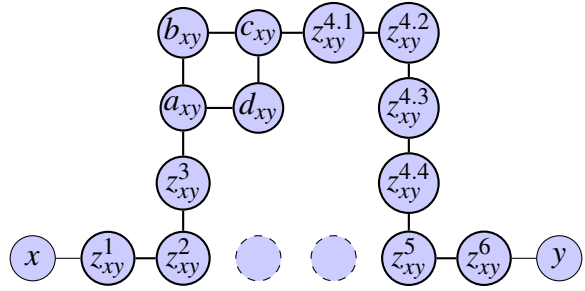

{For every edge $(u,v)\in E$ oriented from $u$ to $v$, we replace, in $H_q$, the first two edges $(u,i), (i,x)$ of $P^{{H_q}}_{uv}$ with $T^2_{ux}$} defined in Figure~\ref{fig:gadget1}, and every other edges $(x,y) \in E^{H_q}$  with $T^1_{xy}$ defined in Figure~\ref{fig:gadget2}.

In the following we use $T_{xy}$ to refer to $T^1_{xy}$ or $T^2_{xy}$.
Note that two gadgets $T_{xy}$ never overlap each other in $F$ and the resulting graph $F$ is a subgrid. Indeed, if $G$ is embedded in a grid $\mathcal G$, $H_q$ is embedded in $Exp(\mathcal G,2q)$ and $F$ is embedded in $Exp(\mathcal G,14q)$.


By construction we have
$|V^F|= |V^{\widetilde H_q}| + 13|E^{\widetilde H_q}| + 3|E|$ and
$|E^F|= 15|E^{\widetilde H_q}| + 3|E|$.

Using Relation~\ref{eq:metricsHtilde}, we deduce:

\begin{equation}\label{eq:metricsF}\begin{array}{lrl}
|V^F|&=&  |V| + 16|E| +28 \sum\limits_{e\in E} k_e\\
|E^F|&=&   18|E| + 30 \sum\limits_{e\in E} k_e
\end{array}\end{equation}

\begin{lem} \label{lem:link_VC_DS}
For any $t\leq |V|$, $\widetilde H_q=(V^{\widetilde H_q}, E^{\widetilde H_q})$ has a vertex cover of size $t$ if and only if $F$ has a dominating set $D$ of size $t+4|E^{\widetilde H_q}|+|E|$ such that, for each edge $(x,y)\in E^{\widetilde H_q}$, we have:
\begin{itemize}
\item at least one vertex of $\{a_{xy},c_{xy}\}$ is in $D$
\item at least one vertex of $\{z^1_{xy},z^{13} _{xy}\}$ if $(x,y)$ is the first edge of a path $P^{\widetilde H_q}_{uv}$ with $(u,v)\in E$ oriented from $u$ to $v$
\item at least one vertex in $\{z^1_{xy},z^{6}_{xy}\}$ in the other cases.
\end{itemize}
\end{lem}

\begin{pf}
For this result, it is convenient to see how $F$ could be constructed from $\widetilde H_q$: for every edge $(u,v)\in E$ oriented from $u$ to $v$, the first edge of $P^{\widetilde H_q}_{uv}$ \-- we denote $E^{\widetilde H_q}_2$ the set of such edges  corresponding to two edges of $P^{ H_q}_{uv}$ \--  is replaced with $T^2_{uv}$. All other edges of $\widetilde H_q$ \-- we denote  $E^{\widetilde H_q}_1 \subset E^{\widetilde H_q}$ their set \-- are replaced with $T^1_{uv}$.  
Note that $|E^{\widetilde H_q}_2|=|E|$.

$\Rightarrow$
Let $U\subset V^{\widetilde H_q}$ be a vertex cover of $\widetilde H_q$ of size $t$. We initialize $D$ with all vertices of $U$, seen as a subset of $V^F$, and complete it in a dominating set of $F$. Then for every $(x,y)\in E^{\widetilde H_q}$, oriented from $x$ to $y$, we have $D\cap \{x,y\} \neq \emptyset$. We then apply one of the two following cases:
\begin{itemize}
    \item[$\bullet$] if $(x,y)\in E^{\widetilde H_q}_2$: If $x\in D$, we add to $D$ the vertices  $z^{3}_{xy},c_{xy},z^{7}_{xy},z^{10}_{xy}$ and $ z^{13}_{xy}$ of $T^2_{xy}$, else if $y\in D$, we add to $D$ the vertices $z^1_{xy},a_{xy},z_{xy}^{6.1},z^{8}_{xy}$ and $z^{11}_{xy}$ of $T^2_{xy}$. In both cases, 5 vertices are added to $D$, and all the vertices of $T^2_{xy}$ are dominated by $D$.
       \item[$\bullet$] if $(x,y)\in E^{\widetilde H_q}_1$: If $x\in D$, we add to $D$ the vertices $z^{3}_{xy},c_{xy},z_{xy}^{4.3}$ and $z^{6}_{xy}$ of $T^1_{xy}$, else if $y\in D$, we add to $D$ the vertices $z^1_{xy},a_{xy},z_{xy}^{4.1}$ and $z_{xy}^{4.4}$ of $T^1_{xy}$. In both cases, 4 vertices are added to $D$, and all the vertices of $T^1_{xy}$ are dominated by $D$.
\end{itemize}
The resulting set $D$ is a dominating set of $F$ of size $t+4|E^{\widetilde H_q}|+|E^{\widetilde H_q}_2|=t+4|E^{\widetilde H_q}|+|E|$ and for each edge $(x,y)\in E^{\widetilde H_q}$, $D$ has at least one vertex in $\{a_{xy},c_{xy}\}$ and one vertex in $\{z^1_{xy},z^{13} _{xy}\}$  (resp.  $\{z^1_{xy},z^{6}_{xy}\}$) if $(x,y)\in E^{\widetilde H_q}_2$(resp. $E^{\widetilde H_q}_1$ ).

$\Leftarrow$
Now suppose we have $D$ a dominating set of $F$. Then for every $(x,y)\in E^{\widetilde H_q}$ oriented from $x$ to $y$, we have:
\begin{itemize}
	\item[$\bullet$] if $(x,y)\in E^{\widetilde H_q}_2$: $D$ includes at least 6~vertices on $T^2_ {xy}$, and 5~vertices on $T^2_ {xy}\setminus \{x,y\}$.
    \item[$\bullet$] if $(x,y)\in E^{\widetilde H_q}_1$: $D$ includes at least 5~vertices on $T^1_ {xy}$, and 4~vertices on $T^1_ {xy}\setminus \{x,y\}$.
\end{itemize}
Then $D$ includes at least $t'+4|E^{\widetilde H_q}|+|E^{\widetilde H_q}_2|=t'+4|E^{\widetilde H_q}|+|E|$ vertices for some integer $t'$. We then perform the following modifications on $D$:
\begin{itemize}
    \item[$\bullet$] for every $(x,y)\in E^{\widetilde H_q}_2$ oriented from $x$ to $y$ : if $x\in D$, we can replace at least 5 vertices of $D\cap T^2_ {xy}$ by $z^{3}_{xy},c_{xy},z^{7}_{xy},z^{10}_{xy}$ and $z^{13}_{xy}$. If $y \in D$, we can replace at least 5 vertices of $D\cap T^2_ {xy}$ by $z^1_{xy},a_{xy},z_{xy}^{6.1},z^{8}_{xy}$ and $z^{11}_{xy}$.
    \item[$\bullet$] for every $(x,y)\in E^{\widetilde H_q}_1$ oriented from $x$ to $y$ : if $x\in D$, we replace at least 4 vertices of $D\cap T^1_ {xy}$ by $z^{3}_{xy},c_{xy},z_{xy}^{4.3}$ and $z^{6}_{xy}$. If $y \in D$, we replace at least 4 vertices of $D\cap T^1_ {xy}$ by $z^1_{xy},a_{xy},z_{xy}^{4.1}$ and $z_{xy}^{4.4}$. If neither $x$ nor $y$ is in $D$, we can induce that $|D\cap T^1_ {xy}|\ge5$. Thus, we replace at least 5 vertices of $D\cap T^1_ {xy}$ by $x,z^{3}_{xy},c_{xy},z_{xy}^{4.3}$ and $z^{6}_{xy}$.
\end{itemize}
Note that none of these modifications increases the size of $D$, and $D$ is still a dominating set of $F$. However, we ensured that $|D\setminus V^{\widetilde H_q}|\ge 4|E|+|E|$, and $|D\cap\{x,y\}|\ge 1, \forall (x,y)\in E^{\widetilde H_q}$. Then $U=D\cap V^{\widetilde H_q}$ is a vertex cover for $\widetilde H_q$ of size at least $t$. This completes the proof. \qed
\end{pf}

\subsection{Relations between {\tt Min P$p$CP} and dominating sets}
\label{subsec:relations_PpCP_DS}

\begin{lem} If $D \subseteq V$ is a dominating set of a triangle-free graph $G=(V,E)$ with degrees $\{2,3\}$, then ${\mathbb E}(D)\le 2$.
\label{lem:relation_DS_Radius}
\end{lem}
\begin{pf}
For any $v\in D$, we recall that $r^s(D,v)=0$ for all scenarios $s$.
For any $v \in V\setminus D$, $v$ is at a distance~1 of a vertex of the dominating set $D$.
Then any neighbor of $v$ is either in $D$ or at a distance~1 of a vertex of $D$ that is not adjacent to $v$, as $G$ is triangle free.
Therefore the evacuation distance of $v$ in any scenario cannot exceed $2$. Thus $r^s(D,v)\le 2$ for all scenarios $s$ and
$${\mathbb E}(D)=\frac{\sum_{s\in V} r^s(D)}{|V|}=\frac{\sum_{s\in V}\max_{v\in V}(r^s(D,v))}{|V|} \le 2$$

%
%
\end{pf}
\qed

\begin{remarque}\label{rem:triangle-free}
 $D$ intersects all articulation components.
\end{remarque}
\begin{lem} For $G=(V,E)$ a 
graph with degrees $\{2,3\}$ and $p<\gamma(G)$, the minimum expected value of the {evacuation radius} over all scenarios of any solution of {\tt Min P$p$CP} is greater than~2.
\label{lem:relation_value_min}
\end{lem}
\begin{pf}
Let $C$ be a solution of {\tt Min P$p$CP} on $G$ for $p<\gamma(G)$. As $C$ cannot be a dominating set, there exists $v\in V$ such that $\{v\}\cup N(v)\cap C = \emptyset$, i.e. $v$ is not adjacent to any vertex of $C$.
For any scenario $s$, the evacuation distance of $v$ will be at least~2 as none of its neighboring vertices is in $C$. Thus $r^s(C,v)\ge 2,  \forall s\in V$, which implies $r^s(C)\ge 2, \forall s\in V$.
In addition, for any vertex $y\in N(v)$, the evacuation distance of $y$ in scenario $y$ is at least~3 since $y$ has an evacuation path that crosses $v$. Since $r^y(C,y)\ge 3$ and $r^y(C)\ge 3$, it follows that ${\mathbb E}(C) > 2$. \qed
\end{pf}

\mar{The following lemma is the counterpart of Lemma~6 in~\cite{DHM18} but requires different arguments.}

\begin{lem}\label{lem:radius}
Let $D$ be a minimum dominating set of $F$ as described in Lemma~\ref{lem:link_VC_DS} and of size $p_t$.
$D$ is a solution of {\tt Min P$p$CP} for $p=p_t$ of value strictly less than~2.
\end{lem}

\begin{pf}
Note that  $|D|=p_t=\tau(G)+4|E^{\widetilde H_q}|+|E^{\widetilde H_q}_2|$ as shown in Lemma~\ref{lem:link_VC_DS}.
Using Remark~\ref{rem:triangle-free} and since $F$ is triangle-free (it is a subgrid), $D$ can then be seen as a feasible solution for {\tt Min P$p$CP} and $p=p_t$ in the graph $F$. We claim the following relation that immediately concludes the proof:

$$r^s(D)=\left\{\begin{array}{rcl}
	1 & {\rm \ if \ }  s\in V^{\widetilde H_q} \subset V^F {\rm \ and \ } s\not\in D\\
	2 & {\rm \  otherwise \ }
\end{array}\right.$$

We recall that every vertex of $\widetilde H_q$ maps a vertex in $F$ by construction, thus we consider $V^{\widetilde H_q}\subset V^F$ in the following.
Since $F$ is triangle-free with no pending vertex, and $D$ is a dominating set, then we have by Lemma~\ref{lem:relation_DS_Radius} $r^s(D,v)\leq 2, \forall s,v\in V^F$. 

Three cases emerge:
\begin{enumerate}
    \item $s \in V^F \setminus V^{\widetilde H_q}$: Denote $(x,y)\in E^{\widetilde H_q}$ such that $s\in T_{xy}$. As $D$ is a minimal dominating set of $F$, $D$ is build as the resulting dominating set described in Lemma~\ref{lem:link_VC_DS}. It follows that there is at least one evacuation distance of length~2 for any scenario $s\in V^F \setminus V^{\widetilde H_q}$, i.e $r^s(D)=2$.
\end{enumerate}

In the following, $s\in V^{\widetilde H_q}$ and we denote by $u_1,\hdots,u_d \in V^{\widetilde H_q}$ the neighbors of $s$ in $\widetilde H_q$.

\begin{enumerate}
    \item[2] $s\in V^{\widetilde H_q} \cap D$: Since $D$ is minimal, $D\cap V^{\widetilde H_q}$ is a minimal vertex cover of $\widetilde H_q$, thus there is at least one neighbor $u \in \{u_1,\hdots,u_d\}$ of $s$ in $\widetilde H_q$ that is not included in $D$. By construction, $z^1_{su}, z^2_{su} \notin D$ and $z^3_{su}\in D$. Then under scenario $s$, the evacuation distance of $z^1_{su}$ is 2, i.e. $r^s(D,{z^1_{su}})=2$. Under scenario $s$, the evacuation distance of any other vertex in $T_{su}$ is less than~2 given that $D$ is a minimal dominating set. For any other neighbor $u'\in \{u_1,\hdots,u_d\}$ of $s$ in $\widetilde H_q$ ($u'\neq u$), we have $|\{z^1_{su'},z^2_{su'},z^3_{su'}\} \cap D|=1$, and $D$ a minimal dominating set on $T_{su'}$, thus the evacuation distance of any vertex in $T_{su'}$ is at most~2. Therefore $r^s(D)=2$.
    \item[3] $s\in V^{\widetilde H_q} \setminus D$: We recall that by definition $D\cap V^{\widetilde H_q}$ is a minimal vertex cover of $\widetilde H_q$, then $\{u_1,\hdots,u_d\} \subset D$. In addition, for any edge $(s,u) \in E^{\widetilde H_q}$ oriented from $s$ to $u$, $D$ includes by construction $z^1_{su}$. Then every neighbor of $s$ in $F$ is included in $D$ by construction.
    Therefore, $r^s(D,s)=1$. Since $D$ is a dominating set in $F$, it remains a dominating set in $F\setminus\{s\}$, which guarantees $r^s(D,v)=1, \forall v \in V^F\setminus\{s\}$. Thus $r^s(D)=1$.
\end{enumerate}{}

So, in all cases except the last one, $r^s(U)=2$, and the proof is complete.\qed \end{pf}

We now are ready to prove the main result of this section.

\subsection{The theorem}
\label{subsec:theorem_complexity}


We will use the following easy lemma proved in~\cite{DHM18}.
\begin{lem} (\cite{DHM18})\\
\label{lem:nopending}
The {\tt Min Vertex Cover} problem is NP-hard in planar graphs with vertices of degree~2 or~3.	
\end{lem}

\begin{theoreme}\label{thm:complexity}
	If P$\neq$ NP, there is no polynomial time approximation for  {\tt Min P$p$CP}  guaranteeing a ratio less than $ \frac{56}{55}$ for subgrids with  vertex degrees~2 or~3, even in the uniform case  (all edge lengths are~1).
\end{theoreme}

\begin{pf}

The proof is by contradiction. Let us suppose there is a polynomial approximation algorithm $\mathtt{A}$ for uniform {\tt Min P$p$CP} which guarantees the approximation ratio~$\rho$ satisfying $1 < \rho < \frac{56}{55}$, on subgrids with vertex degrees~2 or~3 for a parameter $p$.
We will show how to use this algorithm to solve the {\tt Min Vertex Cover} problem on  planar graphs. Lemma~\ref{lem:nopending} gives the contradiction, unless P=NP.

Suppose $\varepsilon>0$ such that $\rho <\frac{56+2\varepsilon}{55+2\varepsilon} <\frac{56}{55}$. Take an integer $q\ge 2$  such that $\varepsilon \ge \frac{17}{q-1}$.

Consider a  planar graph $G=(V,E)$, instance of {\tt Min Vertex Cover}. Consider the graph ${H_q}$ obtained by Transformation~\ref{transformation:1}, as well as $\widetilde H_q=(V^{\widetilde H_q},E^{\widetilde H_q})$ and the vector $\{k_e :e \in E\}$ obtained through ${H_q}$.
In addition, consider the graph $F=(V^F, E^F)$ obtained from ${H_q}$ through Transformation~\ref{transformation:3}.

{Recall that, from Relations~\ref{eq:metricsHtilde} and~\ref{eq:metricsF}, we have $|V^{\widetilde H_q}|=|V|+2\sum\limits_{e\in E} k_e$, $|E^{\widetilde H_q}|=|E|+2\sum\limits_{e\in E} k_e$ and $|V^F|= |V| + 16|E| +28 \sum\limits_{e\in E} k_e$.}
%

%
We also deduce from  Lemma~\ref{lem:link_VC_DS}:

\begin{equation}\label{eq:opt}\begin{array}{lrl}
\gamma(F)&=&\tau(\widetilde H_q)+4|E^{\widetilde H_q}| + |E|  \\
&=&\tau(G)+ 5|E| +9 \sum\limits_{e\in E} k_e
\end{array}\end{equation}

We apply the hypothetical approximation algorithm $\mathtt{A}$ on $F$ for different values of $p$, starting with $p=2$ and augmenting it.
Suppose first we use $p< \gamma(F)$ and the algorithm computes a solution $C$. Then ${\mathbb E}(C)\ge 2$ as proven in Lemma~\ref{lem:relation_value_min}.
Suppose now we set $p= \gamma(F)=\tau(G)+ 5|E| +9 \sum_{e\in E} k_e$.
Given Lemma~\ref{lem:radius}, we obtain the following:
 $$
 {\mathbb E}(C)= \frac{(|V^{\widetilde H_q}|-\tau(\widetilde H_q)) + 2(|V^F|- (|V^{\widetilde H_q}|-\tau(\widetilde H_q)))} {|V^F|} = \frac{2|V^F|-(|V^{\widetilde H_q}|-\tau(\widetilde H_q))} {|V^F|}
 $$

We deduce, using Relations~\ref{eq:tau_G_wideH}, \ref{eq:metricsF} and~\ref{eq:opt}:
$$\begin{array}{rcl}
|V^F|{\mathbb E}(C)&=&|V| + 32|E| + 55\sum\limits_{e\in E} k_e + \tau(G) \\
&<&2|V| + 32|E| + 55\sum\limits_{e\in E} k_e
\end{array}$$

where the last inequality holds because $\tau(G)<|V|$.
So, we have:

$$
{\mathbb E}(C)<\frac{2|V| + 32|E| + 55\sum\limits_{e\in E} k_e}{ |V| + 16|E| +28 \sum\limits_{e\in E} k_e }=2-\frac{\sum\limits_{e\in E} k_e}{ |V| + 16|E| +28 \sum\limits_{e\in E} k_e } $$

Using Equation~\ref{eq:q}, we have   $\sum\limits_{e\in E} k_e \ge (q-1)|E|$. In addition, since $G$ is of degree~2 or~3, we have $|V| \le |E|$. It follows:

$$
{\mathbb E}(C)<2-\frac{\sum\limits_{e\in E} k_e}{17|E| +28\sum\limits_{e\in E} k_e} \le 2 - \frac{1}{28 +  \frac{17|E|}{\sum\limits_{e\in E} k_e}} \le 2 - \frac{1}{28 +   \frac{17}{q-1}}
$$

As $\varepsilon \ge \frac{17}{q-1}$ we get:

$$
{\mathbb E}(C)\le 2 - \frac{1}{28 + \varepsilon} \le \frac{55 + 2\varepsilon}{28+ \varepsilon}
$$

As a consequence,  and since an optimal probabilistic solution $C^\ast$ will satisfy ${\mathbb E}(C^\ast)\leq {\mathbb E}(C)\le \frac{55 + 2\varepsilon}{28+ \varepsilon}$, the approximation algorithm~$\mathtt{A}$ will determine an approximated solution $C$ in $F$ of value:

\begin{equation}\label{eq:lessthan2}
\begin{array}{rcl}
{\mathbb E}(C)&\leq&  \rho \times {\mathbb E}(C^\ast) \\
&\le& \rho \times\frac{55 +2 \varepsilon}{28 + \varepsilon}\\
&<& \frac{56 +2 \varepsilon}{55 + 2\varepsilon} \times \frac{55 +2 \varepsilon}{28 + \varepsilon}\\
&<&2 \\
\end{array}
\end{equation}

Note that, given a solution $C$, computing its probabilistic {radius} can be done in polynomial time.
Indeed, for any $v,s\in V^F$, computing $r^s(C,v)$ can be performed using any minimum path algorithm.
Hence, we can apply successively the approximation algorithm $\mathtt{A}$ on the graph $F$ for increasing values of $p$, starting with $p=2$, until the computed solution $C$ satisfies ${\mathbb E}(C)<2$.
Thanks to Lemma~\ref{lem:relation_value_min} and Equation~\ref{eq:lessthan2}, the algorithm stops for $p= \gamma(F)=\tau(G)+ 5|E| +9 \sum\limits_{e\in E} k_e$.
Using Equations~\ref{eq:opt} we can deduce $\tau(G)= p -5|E| -9 \sum\limits_{e\in E} k_e $.

Since constructing $\widetilde H_q$ and $F$, as well as evaluating ${\mathbb E}(C)$, can be done in polynomial time, and since algorithm~$\mathtt{A}$ will be run less than $|V|$ times, the whole process is polynomial. This is a contradiction if P$\neq$NP, and the proof is complete.
\qed
\end{pf}

\section{Approximation results in the uniform case.}
\label{sec:approximation}


We will show that, in graphs of bounded average degree, there is a polynomial approximation algorithm guaranteeing a constant approximation ratio for the uniform  {\tt Min P$p$CP} (i.e., with all edge lengths equal to~1). Our result is even valid if edge lengths lie into $[l, 2l]$ for a positive $l$.

Our strategy is to show that, under these assumptions, the ratio $\frac{{\mathbb E}(C)}{r(C)}$ is bounded for any  $p$-center $C$ that is feasible for {\tt Min P$p$CP}. In particular, a solution with constant approximation ratio for {\tt Min $p$-Center} \mar{that is feasible for {\tt Min P$p$CP}} has a constant ratio for the latter.

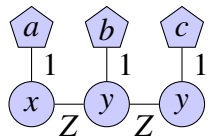
\begin{figure}[hbt]
\centering
 \begin{tikzpicture}
	
	\node[sommet] (v1) at (0,0){$x$} ;
	\node[sommet] (v2) at (1,0){$y$} ;
	\node[sommet] (v3) at (2,0){$y$} ;
	\node[centre] (v4) at (0,1){$a$} ;
	\node[centre] (v5) at (1,1){$b$} ;
	
	\node[centre] (v6) at (2,1){$c$} ;
	
 	\foreach \from/\to in {v1/v2,v2/v3}
  	\draw (\from) edge node[below]{$Z$} (\to) ;
  		\foreach \from/\to in {v1/v4,v2/v5,v3/v6}
  	\draw (\from) edge node[right]{$1$} (\to) ;

 \end{tikzpicture}
\caption{A case where $\frac{{\mathbb E}(C)}{r(C)}=Z+1$. }
\label{fig:L+1}
\end{figure}

Note that in graphs with general lengths we cannot expect the same and thus, another strategy should be taken. Indeed, consider the caterpillar $H$ of Figure~\ref{fig:L+1}  with three internal vertices $x,y,z$ and  edges $(x,y)$ and $(y,z)$ of length $Z$ and three pendent vertices $a,b,c$, respectively linked to $x,y,z$ with edges of length~1.

$\{a,b,c\}$ is the unique feasible solution of the {\tt Min P$p$CP}-instance $(H,3)$. We have\\ $r(\{a,b,c\})=1$. However, for any scenario $s$, $r^s(\{a,b,c\})=Z+1$, which implies\\ ${\mathbb E}(\{a,b,c\})=Z+1$.

Given an edge-weighted graph $G=(V,E,L)$, recall that $\mathcal{C}_p(G)$ denotes the set of feasible solutions of the {\tt Min P$p$CP}-instance $(G,p)$.  From Proposition~\ref{prop:feasibility_condition}, a set $C\subset V$  is in $\mathcal{C}_p(G)$ if and only if $|C|\geq 2$ and $C$ intersects all MACs. For any $p\leq |V|$, we call {\em MAC $p$-center} a $p$-center intersecting all MACs. For $p\geq 2$,  $\mathcal{C}_p(G)$ is the set of MAC $p$-centers.

For any set $C\subset V$ of centers, 
recall that the radius of $C$ is $r(C)=\max\limits_{v\in V} d(v,C)$. Note that  for any scenario $s\in V$, $r^s(C) \ge r(C)$.
We consider the {\em {\tt Min MAC $p$-Center}}  problem of finding a MAC $p$-center of minimum radius. The {\tt Min MAC $p$-Center} problem has a feasible solution for a graph $G$ if and only if $p$ is at least the number of MACs in $G$, i.e., $p\geq |\Upsilon(G)|$.

In what follows, we describe an approximation preserving reduction between {\tt Min P$p$CP} and {\tt Min MAC $p$-Center}  (Subsection~\ref{subsec: approx-reduc}). A polynomial approximation algorithm for the latter leads to a polynomial approximation algorithm for the former with a ratio that depends on the average degree $\overline{deg}(G)=\frac{2|E|}{|V|}$ of $G$. More precisely, the reduction is even the identity and we analyze how good {for the problem {\tt Min P$p$CP}} an approximated MAC $p$-center can be. Then, in Subsection~\ref{subsec: approx-pcenter}, we show that {\tt Min MAC $p$-Center}  can be approximated within the ratio 2, which leads to a $(4\overline{deg}(G) +2)$-approximation for the uniform {\tt Min P$p$CP} (all edges are of length~1). {Actually, the result  still holds if all edge-lengths lie in the interval} $[\ell, 2\ell]$ for any positive $\ell$. 

\subsection{A polynomial approximation preserving reduction}\label{subsec: approx-reduc}

We directly establish the following proposition for general edge lengths. We will denote respectively $\ell_M$ and $\ell_m$ the maximum and minimum edge lengths. 

\begin{proposition}
\label{prop:approx_1}
On an edge weighted graph with lengths in $[\ell_m, \ell_M]$,  $ \forall C\in \mathcal{C}_p(G)$, we have:
$${\mathbb E}(C) \le (2\overline{deg}(G)+1)r(C)+(\ell_M-2\ell_m)\overline{deg}(G)$$
\end{proposition}
\begin{pf}
Let us consider any scenario $s\in V$ of degree $deg(s)$ and number $1,2,\hdots, deg(s)$ the edges incident to $s$.
We claim that $r^s(C) \le (2deg(s)+1)r(C)$.

Consider indeed $x\in V$ such that $r^s(C,x) = r^s(C) \ge r(C)$.
If $r^s(C,x)=r(C)$, then the claim is satisfied.
Let us assume $r^s(C,x)>r(C)$.
We consider two cases.

{\underline{Case 1}}: $x\neq s$.
$r^s(C,x)$ is the  length  of a path $\mu=[x_0,x_1,\hdots, x_k]$, where $x_0=x$, $x_k\in C$ and $\mu$ is a minimum path in $G^s$.



Since $d^s(x,x_k)>r(C)$, we can define $i=\max\{j\in \{0, \ldots k-1\},d^s(x_j,x_k)>r(C)\}$.


Then all vertices $x_j, j \in \{ 0, \hdots, i\}$ are, in $G$, at distance at most $r(C)$ from $s$.
Indeed, the path $x_j,\hdots, x_k$ is a minimum path of length greater than $r(C)$ in $G^s$.
So, in $G$, the evacuation path of vertices $x_j, j \in \{ 0, \hdots, i\}$ passes through $s$.



Figure~\ref{fig:prop_approx_1} illustrates the distance  relation between $x$, $s$ and $x_k$ in the case $x\neq s$.
In the figure, no shelter is located on $s$, but the reasoning is the same if there is one.

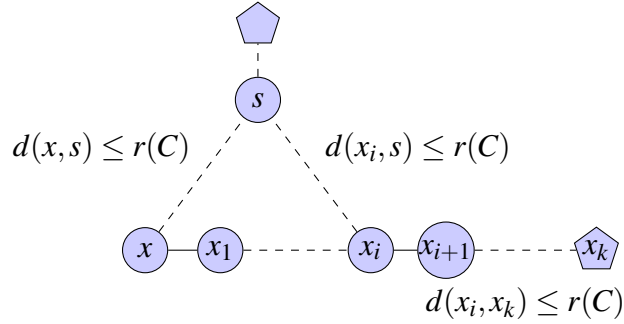
\begin{figure}[hbt]
\centering
 \begin{tikzpicture}
	
	\node[sommet] (v1) at (0,0){$x$} ;
	\node[sommet] (v2) at (1,0){$x_1$} ;
	\node[sommet] (v3) at (3,0){$x_i$} ;
	\node[sommet] (v4) at (4,0){$x_{i+1}$} ;
	\node[centre] (v5) at (6,0){$x_k$} ;
	
	\node[sommet] (v6) at (1.5,2){$s$} ;
	
	\node[centre] (v7) at (1.5,3){ } ;
	
 	\foreach \from/\to in {v1/v2,v3/v4}
  	\draw (\from) -- (\to) ;
  	
  	\draw [dashed] (v1) edge node[above left]{$d(x,s)\le r(C)$} (v6)  ;
  	\draw [dashed] (v2) edge (v3) ;
  	\draw [dashed] (v4) edge node[below=10pt]{$d(x_i,x_k)\le r(C)$}(v5) ;
  	\draw [dashed] (v3) edge node[above right]{$d(x_i,s)\le r(C)$} (v6) ;
  	\draw [dashed] (v6) edge (v7) ;
  	

 \end{tikzpicture}
\caption{Distance relations between vertices $x$, $s$ and $x_k$ used for Proposition~\ref{prop:approx_1}}
\label{fig:prop_approx_1}
\end{figure}

In $G$, for $j\in \{0,\hdots,i\}$, we consider a minimum path from $x_j$ to $s$, of value at most $r(C)$.
We assign to $x_j$ a color in $N_1, \hdots, N_{deg(s)}$ depending on the last edge of the minimum path we have fixed for $x_j$:
$x_j$ is of color $N_t$ if the related minimum path between $x_j$ and $s$ terminates with the $t^{th}$ edge incident to $s$.

Note that the distance in $G_s$ between two vertices of the same color is at most $2r(C)-2\ell_m$.
Indeed considering, in $G$, two minimum paths from these vertices to $s$ and sharing the last edge, we deduce a walk avoiding $s$ between them  of total length at most $2r(C)-2\ell_m$.
This walk includes a path in $G^s$ of length at most $2r(C)-2\ell_m$ between these two vertices.

This allows us to derive an upper bound of $d^s(x,x_i)$.
Suppose $x$ is of color $N_{i_1}$ and consider the last vertex $x_j$ of color $N_{i_1}$  along the path $\mu$; we have $d^s(x,x_j)\leq 2r(C)-2\ell_m$. Then, if $j<i$, the vertex $x_{j+1}$ is of color $N_{i_2}$ and $d^s(x_j, x_{j+1})\leq \ell_M$. Using the same reasoning for all non-empty colors gives $d^s(x,x_i)\leq  deg(s)(2r(C)-2\ell_m) + (deg(s)-1)\ell_M$.

Taking into account the edge $x_ix_{i+1}$ and the fact that $d^s(x_{i+1},x_k)\le r(C)$ we have:

\begin{equation}\label{eq:deg1}
r^s(C) \le    (2deg(s)+1) r(C) +deg(s)(\ell_M-2\ell_m)
\end{equation}

{\underline{Case 2}}: $x= s$
Similarly, $r^s(C,s)$ is the  length of a path $\mu=[x_0,x_1,\hdots, x_k]$, where  $x_0=s$, $x_k\in C$ and $[x_1,\hdots, x_k]$ is a minimum path in $G^s$. We define  $i$
as in the previous case and use the same argument:   $x_1$ is   color $N_{i_1}$ and we define $x_j$ as previously. The only difference is that for any vertex the fixed minimum path from  $x_j$ to $s$ passes through $x_1$ and consequently $d^s(x_1,x_j)\le r(C)-\ell_m$. For the other colors, the same bound as previously holds.  We then get a better bound:
\begin{equation}
\begin{split}
r^s(C) \le &  r(C)-\ell_m + (deg(s)-1)(2r(C)-2\ell_m) + deg(s)\ell_M + r(C)\\
\le & 2deg(s) r(C) +deg(s)(\ell_M-2\ell_m) +\ell_m
\end{split}
\end{equation}

This bound is better than in Equation~\ref{eq:deg1} since $
\ell_m\leq r(C)$.
So, in all cases we have $r^s(C) \le (2deg(s)+1)r(C)$.
We deduce, by taking the average value, ${\mathbb E}(C) =  \frac{1}{|V|}\sum_{s\in V}  r^s(C) \le (2 \overline{deg}(G)+1)r(C)+(\ell_M-2\ell_m)\overline{deg}(G)$ which concludes the proof.
\qed \end{pf}

On a tree, the analysis can be improved:
\begin{proposition}\label{prop:reductree}
On a tree with edge lengths in $[\ell_m, \ell_M]$,  $ \forall C\in \mathcal{C}_p(G)$, we have:
$${\mathbb E}(C)\leq 3r(C)+\ell_M-2\ell_m$$
\end{proposition}
\begin{pf}
Consider, for a scenario $s$, and a vertex $x, r^s(C,x)=r^s(C)$, the same analysis as in the proof of Proposition~\ref{prop:approx_1}.  Since there is no cycle, all vertices $x, \ldots, x_i$ are of the same color. Equation~\ref{eq:deg1} becomes
$$r^s(C)\le 3r(C)+ \ell_M-2\ell_m$$ which concludes the proof.
\qed \end{pf}

\begin{remarque}
In~\cite{DHM18}, we have shown that, on paths with all edge-weight~1, there is an optimal solution $C^*$ of  {\tt Min MAC $p$-Center} such that ${\mathbb E}(C^*)=r(C^*)$.
\end{remarque}

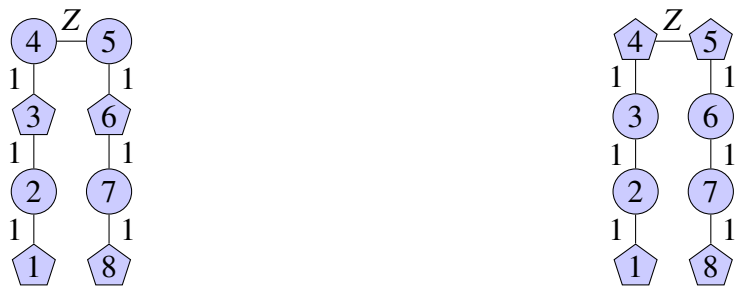
\begin{figure}[hbt]
\centering
 \begin{tikzpicture}
	
	\node[centre] (v1) at (-4.5,0){1} ;
	\node[sommet] (v2) at (-4.5,1){2} ;
	\node[centre] (v3) at (-4.5,2){3} ;
	\node[sommet] (v4) at (-4.5,3){4} ;
	\node[sommet] (v5) at (-3.5,3){5} ;
	\node[centre] (v6) at (-3.5,2){6} ;
	\node[sommet] (v7) at (-3.5,1){7} ;
	\node[centre] (v8) at (-3.5,0){8} ;
	
 	\foreach \from/\to in {v1/v2,v2/v3,v3/v4}
  	\draw (\from) edge node[left]{1} (\to) ;
  	\foreach \from/\to in {v5/v6,v6/v7,v7/v8}
  	\draw (\from) edge node[right]{1} (\to) ;
	\draw (v4) edge node [above]{$Z$} (v5);

  \node[align=left,text width=7cm] at (-4,-1.5) {\footnotesize{$\{1,3,6,8\}$ is an optimal MAC 4-center\\
$r(\{1,3,6,8\})=1,  {\mathbb E}(\{1,3,6,8\})=\frac{Z}{2}+1$}};
	
	\node[centre] (v11) at (3.5,0){1} ;
	\node[sommet] (v12) at (3.5,1){2} ;
	\node[sommet] (v13) at (3.5,2){3} ;
	\node[centre] (v14) at (3.5,3){4} ;
	\node[centre] (v15) at (4.5,3){5} ;
	\node[sommet] (v16) at (4.5,2){6} ;
	\node[sommet] (v17) at (4.5,1){7} ;
	\node[centre] (v18) at (4.5,0){8} ;
	
 	\foreach \from/\to in {v11/v12,v12/v13,v13/v14}
  	\draw (\from) edge node[left]{1} (\to) ;
  	\foreach \from/\to in {v15/v16,v16/v17,v17/v18}
  	\draw (\from) edge node[right]{1} (\to) ;
	\draw (v14) edge node [above]{$Z$} (v15);
\node[align=left,text width=7cm] at (4,-1.5) {\footnotesize{$\{1,4,5,8\}$ is {\tt Min P$p$CP}-optimal for $p=4$;\\
$r(\{1,4,5,8\})=2,  {\mathbb E}(\{1,4,5,8\})=2$}};

 \end{tikzpicture}
\caption{With general weights, an optimal MAC $p$-center can be a very bad {\tt Min P$p$CP} solution.}

\label{fig:P8}
\end{figure}

As noticed in the following example in Figure~\ref{fig:P8}, with general weight system the situation may be totally different. In this example the graph is a path on 8 vertices with only one edge of weight $Z>1$ and all other edges of weight~1 and $p=4$. There is a unique optimal MAC 4-center and, for large values of~$Z$, its value is very bad compared to an optimal {\tt Min P$p$CP} solution.

\begin{proposition}\label{prop:reduction}
Suppose a class of edge-weighted graphs $G=(V,E,L)$ with $\ell_M\leq 2\ell_m$ for which  {\tt Min MAC $p$-Center} can be approximated with $\rho(G)$. \\
Then, {\tt Min P$p$CP} can be approximated with $(2\overline{deg}(G)+1)\rho(G)$ on the same class.
\end{proposition}
\begin{pf}
Given a graph $G$ in the class, we build a $p$-center $C$ in $\mathcal{C}_p(G)$, if it exists, of value at most $\rho(G)r^*(G)$, where $r^*(G)$ denotes the optimal radius of a MAC $p$-center in $G$.
Using Proposition~\ref{prop:approx_1} and $\ell_M\leq 2\ell_m$, we have ${\mathbb E}(C)\le (2\overline{deg}(G)+1)r(C)\le (2\overline{deg}(G)+1)\rho(G)r^*(G)$.

Now if $C^*$ is an optimum solution for {\tt Min P$p$CP}, we have ${\mathbb E}(C^*) \ge r(C^*) \ge r^*(G)$.
This concludes the proof.
\qed \end{pf}

\subsection{Constant approximation algorithms}\label{subsec: approx-pcenter}

The main objective of this subsection is to derive constant approximation results for {\tt Min MAC $p$-Center} using Proposition~\ref{prop:reduction}.
The following easy remark on trees will allow to immediately deduce a first   result on trees.

\begin{proposition}
{\tt Min MAC $p$-Center} is polynomial on trees with general lengths.
\end{proposition}
\begin{pf}
Given a tree $\mathcal{T}$, for any distance $d$ we consider the tree $\mathcal{T}_d$ obtained from $\mathcal{T}$ by gluing  to each pending vertex $v$ a path of length $d$. Then, $\mathcal{T}$ has a MAC $p$-center of radius~$d$ if and only if  $\mathcal{T}_d$ has a $p$-center of radius~$d$. The result immediately follows from the fact that $p$-Center is polynomial on trees. 
\end{pf}

Using Proposition~\ref{prop:reductree} and the analysis of
Proposition~\ref{prop:reduction}, we get:
\begin{corollary}\label{cor:tree}
There is a polynomial algorithm for {\tt Min P$p$CP} guaranteeing the ratio~3 on trees with all edge values~1.
\end{corollary}

\begin{remarque}
 Note however that we leave open the problem of whether {\tt Min P$p$CP} is NP-hard or polynomial on trees.
\end{remarque}

In the reminder of this subsection we devise a 2-approximation polynomial-time  algorithm for {\tt Min MAC $p$-Center} in order to deduce an approximation algorithm for {\tt Min P$p$CP} using Proposition~\ref{prop:reduction}.


\mar{To properly explain the ideas of Algorithm~\ref{algo: mac_pcenter},} we will need another $p$-center problem called
{\em {\tt Min Partial $p$-Center}} that was introduced in~\cite{Partial-p-center}.  Given a graph $G=(V,E)$ and a set of vertices $U\subset V$, {\tt Min Partial $p$-Center}  is to minimize the {\em partial radius} $r(C,U)$ of a $p$-center $C$, where $r(C,U)=\max\limits_{x\in U}d(x,C)$.  The underlying logic is that only vertices in $U$ need to be close to a center. However, centers can be any vertex in $G$ and distances are computed in $G$ (within our terminology, it means that {the evacuation} paths toward a shelter are not required to stay in $U$). 

\mar{The idea of Algorithm~\ref{algo: mac_pcenter} is to reduce {\tt Min MAC $p$-Center} to {\tt Min Partial $p$-Center}  through a pre-processing that allocates some centers to MACs. Then, the solution is completed using {\tt Min Partial $p$-Center}. As we will see, {\tt Min Partial $p$-Center} can be approximated by generalizing the 2-approximation algorithm for {\tt Min $p$-Center}  in~\cite{HS85} or using the general method in~\cite{HS86}. However, since it cannot be directly deduced from existing results, we  will give a direct proof through few claims. }

\subsubsection{Lower bound on the approximation ratios}

Note that, if $U=V$, then $r(C,V)=r(C)$ and {\tt Min Partial $p$-Center} is just the usual {\tt Min $p$-Center} problem. So, {\tt Min $p$-Center} is a particular case or {\tt Min Partial $p$-Center}.  In  particular, {\tt Min Partial $p$-Center}  is \mar{NP-hard and} not approximable within $2-\varepsilon$ for any $\varepsilon >0$, unless P=NP by using the same hardness result for {\tt Min $p$-Center}  proved in~\cite{Hard-Bottleneck}. Note that this hardness result for {\tt Min $p$-Center}, directly obtained from the NP-hardness of {\tt Min Dominating Set},
 holds in the uniform case (all edges have the length~1). Since {\tt Min Dominating Set} remains NP-hard in planar bipartite graphs of degree~3, {\tt Min $p$-Center}, and by consequence {\tt Min Partial $p$-Center}, are not approximable within $2-\varepsilon$ for any $\varepsilon >0$ in planar bipartite graphs of degree~3 with all edge lengths~1, unless P=NP.

Note that the argument used for {\tt Min Partial $p$-Center} cannot be easily adapted to  {\tt Min MAC $p$-Center} since this latter problem is not an immediate  generalization of {\tt Min $p$-Center}. However, for any edge weighted graph $G=(V,E,L)$, instance of  $p$-Center, the instance is equivalent to the instance $(K,\tilde L)$, where $K$ is the complete graph over $V$ and $\tilde L$ denotes the minimum path distance, i.e., $\forall i,j\in V, \tilde \ell_{ij}=d(i,j)$, where the distance $d$ is the distance in $G$. Both instances $G$ and $K$ have the same feasible solutions with the same values and thus, the same optimal solutions. To guarantee finite edge lengths in $K$, we just consider $G$ is connected. Since $K$ is 2-connected as soon as $|V|\geq 2$,  {\tt Min MAC $p$-Center}  is equivalent to {\tt Min $p$-Center}  on $K$. Since the hardness result for {\tt Min $p$-Center}  still holds in connected graphs, {\tt Min MAC $p$-Center} is   not approximable within $2-\varepsilon$ for any $\varepsilon >0$, unless P=NP. We can even easily show that this hardness results already holds for the uniform case where all edge lengths are~1. To this aim, we just need to show that {\tt Min Dominating Set}  is NP-hard in 2-connected graphs. Given a graph $G=(V,E)$  instance of {\tt Min Dominating Set}, we construct $G'$ from $G$ as follows: for every articulation point $a$ of $G$, create a twin vertex $a'$ linked to $a$ and to all neighbors of $a$. $G'$ is 2-connected and the {\tt Min Dominating Set} problems in $G$ and  $G'$ are equivalent. Now, a set of $p$ vertices in $G'$ is a dominating set if and only if its radius is~1 and else, the minimum radius of a $p$-center is at least~2. It immediately implies:

\begin{proposition}\label{prop:lb2mac}
{\tt Min MAC $p$-Center} in graphs with edge lengths all equal to~1 is   not approximable within $2-\varepsilon$ for any $\varepsilon >0$, unless P=NP.
\end{proposition}

\subsubsection{Approximation algorithms}


Consider an instance $(G,U)$ of {\tt Min Partial $p$-Center}, where $G=(V,E,L)$ is a graph with positive lengths on edges and $U\subset V$. We denote $n=|V|$. We can compute $K=(V,\tilde E, \tilde L)$ in $O(n^3)$. We denote  $SL=\{d(x,y), x,y\in V\}$ the set of edge lengths in $K$  (note that $|SL|\leq n^2$) and for any $d\in SL$, $K_d=(V,E_d)$ is the partial graph of $K$ where $E_d$ is the set of edges of length at most $d$. Note that for any $p$-center, its radius is in $SL$.


A  $p$-center of partial radius $d$ in $(G,U)$ can be seen as a {\em partial dominating set} of $(K_d,U)$, where a partial dominating set $X$ is a set of vertices such that every vertex in $U$ has at least one neighbor in $X$. If $A_d$ is the adjacency matrix of $K_d$ with additional 1s on the diagonal (alternatively $A_d$ is the adjacency matrix of $K_d$ with additional loops on each vertex), we denote $A_{d,U}$  the sub-matrix of $A_d$ corresponding to rows in $U$ (it has $|U|$ rows and $|V|$ columns). The problem of finding a minimum partial dominating set can the formulated by the following mathematical program $PDS(G,U,d)$, where the 1s on the diagonal represent the fact that a vertex dominates itself \mar{and the notation $\mathbf{1}_d$ for an integer $d$ denotes a vertical vector of dimension $d$ with only 1-entries}:

$PDS(G,U,d):  \left\{ \begin{array}{rcl} \min & \langle \mathbf{1}_{|V|},x\rangle\\
&A_{d,U}x\geq \mathbf{1}_{|U|}\\
&x\in\{0,1\}^{|V|}
\end{array}\right.
$

We then consider the mathematical program $SIS(G,U,d)$ that corresponds to finding a maximum {\em strong independent set} of $K_d$ contained in $U$, where  a strong independent set $S\subset V$  is an independent set (every two vertices in $S$ are not adjacent) such that every vertex in $V\setminus S$ has at most one neighbor in $S$.

$SIS(G,U,d): \left\{ \begin{array}{rcl} \max & \langle \mathbf{1}_{|U|},y\rangle\\
&A_{d,U}^\intercal y\leq \mathbf{1}_{|V|}\\
&y\in\{0,1\}^{|U|}
\end{array}\right.
$

\begin{monclaim}\label{claim-sis1}
The cardinality of any strong independent set of $K_d$ contained in $U$ is not more than the cardinality of any partial dominating set of $(K_d,U)$.
\end{monclaim}

\begin{pf}
The relaxations of mathematical programs $PDS(G,U,d)$ and $SIS(G,U,d)$, replacing the binary conditions with non negative conditions, are dual linear programming problems. The result is an immediate consequence of  the weak duality theorem.
\qed\end{pf}

Let $d_{max}=\max(SL)$.
We denote $K_{2d,U}$ the graph $K_{\min(2d, d_{max})}[U]$.

\begin{monclaim}\label{claim-sis2}
For a given distance $d\in SL$, let $S_d$ be a maximal independent set  of  $K_{2d,U}$.
$S_d$ is a partial $|S_d|$-center in $(G,U)$ of partial radius $r(S_d,U)\leq 2d$.
\end{monclaim}

\begin{pf}
Consider any vertex $u\in U\setminus S_d$. Since $S_d$ is maximal,
 $S_d\cup \{u\}$ is not independent in $K_{2d,U}$, which  means $d(u,S_d)\leq 2d$ and the claim is proved.
\qed\end{pf}

\begin{monclaim}\label{claim-sis3}
Any independent set $S$ of $K_{2d,U}$ is a strong independent set of $K_d$ contained in $U$.
\end{monclaim}

\begin{pf}
By definition, $S\subset U$. Since $S$ is independent in $K_{2d,U}$,  it is independent in $K_{d,U}$, a partial graph of $K_{2d,U}$. So, it is an independent set of $K_d$. The result then follows by contrapositive: if there is  a vertex $u\in V\setminus S$
 adjacent, in $K_d$, to two vertices of $S$, then these two vertices would be at distance at most $2d$, so would be adjacent in $K_{2d,U}$.
\qed\end{pf}

{Claims~\ref{claim-sis1}, \ref{claim-sis2} and \ref{claim-sis3} immediately allow to derive an approximation algorithm for {\tt Min Partial $p$-Center}.} \mar{Even if this result is not strictly used for Theorem~\ref{th: approx-mac_pcenter}, it is worth to mention and it helps understanding the main ideas of Algorithm~\ref{algo: mac_pcenter}.}

\begin{proposition}\label{th: approx-partial_pcenter}
{\tt Min Partial $p$-Center}  is polynomially 2-approximable and this is the best possible constant ratio.
\end{proposition}

\begin{pf}
We already noted that 2 is a lower bound for any constant approximation ratio of {\tt Min Partial $p$-Center}. So, we only need to prove that this bound can be guaranteed.

For a given instance $(G,U)$, we can compute $SL$ and all distance $d(i,j), i,j\in V$ in $O(n^3)$. Then, for any $d\in SL$, we can compute a maximal  independent set $S_d$ of $K_{2d,U}$ and then select $S_{\tilde d}$, where
$\tilde d\in\argmin\limits_{d\in SL, |S_d|\leq p}(r(S_d))$. In other words, $S_{\tilde d}$ is of minimum value among all $S_d$s of cardinality at most $p$. Denote $r^*_U$ the minimum  partial radius of a $p$-center in $(G,U)$.  $r^*_U\in SL$.  Using Claim~\ref{claim-sis3} and Claim~\ref{claim-sis1}, $|S_{r^*_U}|\leq p$ and thus, $\tilde d$ exists and $r(S_{\tilde d})\leq r(S_{r^*_U})$. Using Claim~\ref{claim-sis2}, we deduce $r(S_{r^*_U})\leq 2r^*_U$, which completes the proof.
\qed \end{pf}

Note that, using a binary search on the same model as the 2-approximation algorithm for {\tt Min $p$-Center} proposed in~\cite{HS85}, we can  design a 2-approxi\-ma\-tion algorithm of complexity $O(n^2\log n)$ as soon as all distances between two vertices in $G$ are computed.

{We use similar ideas and the same claims to derive a polynomial 2-approximation algorithm for {\tt Min MAC $p$-Center} (Algorithm~\ref{algo: mac_pcenter}).}

\begin{algorithm}[hbtp]
\floatname{algorithm}{Algorithm}
\caption{2-approximation for {\tt Min MAC $p$-Center} and {\tt Min P$p$CP}.}
\begin{algorithmic}[1]\label{algo: mac_pcenter}
\REQUIRE{Edge weighted graph $G=(V,E,L)$ (lengths are non negative) and $p\geq 2$.}
\ENSURE{Outputs $C$, a MAC $p$-center if it exists.}
\STATE {\bf Begin}
\STATE{Compute $A_1, \ldots A_k$, and $a_1, \ldots a_k$ }
\IF{$k>p$}
\STATE{No-solution output}
\ELSE{}
\STATE{Compute $SL$ and all distances $d(i,j), i,j\in V$}
\STATE{$\widetilde{SL}\leftarrow \emptyset$}
\FOR{$d\in SL$}
\STATE{Compute $I_d^-$ and $I_d^+$}
\STATE{$C_d\leftarrow \emptyset$}\label{line:Cd start}
\FOR{$i\in I_d^-$}
\STATE{Select $x\in A_i$ }
\STATE{$C_d \leftarrow C_d\cup \{x\}$}
\ENDFOR \label{line:endfor-}
\STATE{$V'_d\leftarrow \{v\in V, d(v,\{a_i, i\in I_d^-\})>d\}$}\label{line:V'd}
\STATE{$S_d\leftarrow \emptyset$}\label{line:for+}
\FOR{$i\in I_d^+$}
\STATE{Select $y\in \argmax\limits_{x\in A_i}d(x,a_i)$ }\label{line:Id+}
\STATE{$S_d \leftarrow S_d\cup \{y\}$}
\ENDFOR \label{line:endfor+}
\WHILE{$\exists v\in V'_d, d(v,S_d)>2d$}\label{line:while}
\STATE{$S_d\leftarrow S_d\cup\{v\}$}
\ENDWHILE \label{line:endwhile}
\IF{$|S_d|\leq p-|I_d^-|$}
\STATE{$\widetilde {SL}\leftarrow \widetilde {SL}\cup\{d\}$}\label{line:Dtilde}
\STATE{$C_d \leftarrow C_d\cup S_d$}
\ENDIF
\ENDFOR \label{line:endfor}
\STATE{Let $\tilde d \in \argmin\limits_{d\in \widetilde{SL}}(r(C_d))$}\label{line:tilded}
\STATE{$C\leftarrow C_{\tilde d}$}
\RETURN{$C$}
\ENDIF
\STATE {\bf End}

\end{algorithmic}

\end{algorithm}

To simplify the description of Algorithm~\ref{algo: mac_pcenter}, we introduce some notations used in the description of the algorithm.
Given the instance $G=(V,E,L)$, we denote by $k$ the number of MACs of $G$. These MACs are denoted  $A_1, \ldots A_k$ and the related articulation points are called $a_1, \ldots a_k$  (we may have $a_i=a_j, i\neq j$). As previously $SL=\{d(i,j), i,j\in V,\}$; for any $d\in SL$, we partition $I=\{1, \ldots, k\}$ into $I=I_d^-\sqcup I_d^+$  ($\sqcup$ denotes the disjoint union), where $I_d^-=\{i\in I, \max\limits_{x\in A_i}d(x,a_i)\leq d\}$ and $I_d^+=\{i\in I, \max\limits_{x\in A_i}d(x,a_i)> d\}$. {MACs $A_i, i\in I_d^-$ are seen as small MACs relative to $d$, while MACs $A_i, i\in I_d^+$ are seen as large ones.}
``No-solution output'' is any output we  use to indicate that the problem has no feasible solution.

{The idea of the Algorithm is as follows:
\begin{enumerate}
\item If the number of MAC is more than $p$, then there is obviously no solution.
    \item Else, for every distance $d\in SL$, Algorithm~\ref{algo: mac_pcenter} tries to compute a MAC $p$-center  $C_d$ of radius at most $2d$; only  feasible MAC $p$-centers obtained through this process will be kept and $\widetilde {SL}$ is the set of distances $d$ for which it will occur;
    \item $C_d$ is built as follows: 
    \begin{enumerate}
       \item 
   The algorithm selects one center per small MAC $A_i, i\in I_d^-$;
    \item For each $i\in I_d^-$, all vertices at distance at most $d$ from $a_i$ are allocated to the related center (by definition of $I_d^-$, this includes in particular all vertices of $A_i$).
    \item $V'_d$ is the set of uncovered vertices. If  possible, the algorithm completes  $C_d$ with a partial $(p-|I_d^-|)$-center of $(G\setminus \bigcup\limits_{i\in I_{d}^-}A_i ,V'_{d})$ of partial radius at most $2d$. To this aim, it uses the same ideas as in Proposition~\ref{th: approx-partial_pcenter}: it constructs a maximal independent set $S_d$ of $K_{2d,V'_d}$, but to ensure it intersects all $A_i$s, $i\in I_{d}^+$, it initializes it by choosing one vertex in each of these components. If $|S_d|\leq p-|I_d^-|$, then $d\in \widetilde {SL}$; 
    \end{enumerate}
    \item The best solution $C_{\tilde d}, d\in \widetilde {SL}$ is selected as an approximated solution for {\tt Min MAC $p$-Center}.
\end{enumerate}
}

\begin{theoreme}\label{th: approx-mac_pcenter}
Algorithm~\ref{algo: mac_pcenter} is a polynomial 2-approximation algorithm for {\tt Min MAC $p$-Cen\-ter} and this is the best possible constant ratio.
\end{theoreme}

\begin{pf}
We already noted that 2 is a lower bound for constant approximation ratios. So, we only need to prove that this bound can be guaranteed.

Assume that $k\leq p$; then the instance of {\tt Min MAC $p$-Center} has feasible solutions and thus, also an optimal solution.



Fix a distance $d\in SL$. Note first that, by definition of $I_d^-$ and  $I_d^+$, $V'_d$ computed at line~\ref{line:V'd} satisfies $V'_d\subset V\setminus \bigcup\limits_{i\in I_d^-}A_i$ and $\forall i\in I_d^+$, $A_i\cap V'_d\neq\emptyset$.
Then, the algorithm computes the set $S_d$  from Lines~\ref{line:for+} to Line~\ref{line:endwhile}.

\begin{monclaim}\label{claim:mac1}
$\forall d\in SL$, $S_d$ is a maximal independent set in $K_{2d,V'_d}$ that intersects all $A_i$s, $i\in I_d^+$.
\end{monclaim}

\begin{pf}
The algorithm initializes $S_d$ by selecting,  in each MAC $A_i, i\in I_d^+$, a vertex at maximum distance from $a_i$. This ensures that, at Line~\ref{line:endfor+}, $S_d$ includes one element per MAC $A_i, i\in I_d^+$ and is an independent set (possibly empty) in $K_{2d,V'_d}$. Indeed, if $y_i, y_j$ are respectively selected at Line~\ref{line:Id+} for $i,j\in I_d^+, i\neq j$, then any path between them passes through $a_i$ and $a_j$ (we may have $a_i=a_j$) and is of length greater than $2d$. As a consequence, $S_d$  is a maximal independent set in $K_{2d,V_d'}$.
\qed \end{pf}

$\widetilde {SL}$, computed by the algorithm (Lines~\ref{line:Dtilde}), is the set of distances $d$  such that $S_d$ is of size at most $p-|I_d^-|$.
Consider now an optimal MAC $p$-center, $C_{MAC}^*$, of radius $d^*$.

\begin{monclaim}\label{claim:mac2}
$d^*\in \widetilde {SL}$
 \end{monclaim}

\begin{pf}
Since $C_{MAC}^*$ has at least one center per MAC, $C_{MAC}^*$ has at most $\left(p-|I_{d^*}^-|\right)$ centers in $V\setminus \bigcup\limits_{i\in I_{d^*}^-}A_i$. In addition, vertices in $V'_{d^*}$ cannot be associated with (i.e., evacuated to) centers in $\bigcup\limits_{i\in I_{d^*}^-}A_i$ since these centers are at distance more than $d^*$. This means that $C_{MAC}^* \cap (V\setminus \bigcup\limits_{i\in I_{d^*}^-}A_i)$ is a $(p-|I_{d^*}^-|)$-center of partial radius at most
$d^*$ in $(G\setminus \bigcup\limits_{i\in I_{d^*}^-}A_i ,V'_{d^*})$.

As a consequence $C_{MAC}^* \cap (V\setminus \bigcup\limits_{i\in I_{d^*}^-}A_i)$ is a partial dominating set in $(K_{d^*}, V'_{d^*})$. Using Claims~\ref{claim-sis3} and~\ref{claim-sis1}, we get $|S_{d^*}|\leq |C_{MAC}^* \cap (V\setminus \bigcup\limits_{i\in I_{d^*}^-}A_i)|\leq p-|I_{d^*}^-|$, which means $d^*\in \widetilde {SL}$.
\qed \end{pf}

Claim~\ref{claim:mac2} ensures in particular that $\widetilde {SL}\neq\emptyset$ and consequently $\tilde d$ computed at Line~\ref{line:tilded} is well defined. Since $d^*$ and $\tilde d$ are both in $\widetilde {SL}$, the algorithm computes both sets $C_{d^*}$ and $C_{\tilde d}$ by selecting one vertex per $A_i, i\in I_{d^*}^-$ and one vertex per $A_i, i\in I_{\tilde d}^-$, respectively (from Line~\ref{line:Cd start} to Line~\ref{line:endfor-}) and completing with $S_{d^*}$ and $S_{\tilde d}$, respectively. 
Using Claim~\ref{claim:mac1}, this ensures that both $C_{d^*}$ and $C_{\tilde d}$ are MAC $p$-centers.

Finally, $C_{\tilde d}$ is selected as approximated  solution and Line~\ref{line:tilded} ensures

  \begin{equation}\label{eq:tilde-star}
r(C_{\tilde d})\leq r(C_{d^*})
\end{equation}

We complete the proof by showing the following claim.

\begin{monclaim}\label{claim:mac3}
$r(C_{d^*})\leq 2d^*$.
\end{monclaim}

\begin{pf}
 Consider first a vertex  $v\in V'_{d^*}$ and {use the same argument as in the proof of Proposition~\ref{th: approx-partial_pcenter}.} We have $d(v,C_{d^*})\leq d(v,C_{d^*}\setminus \bigcup\limits_{i\in I_{d^*}^-}A_i )\leq r(S_{d^*},V'_{d^*})$.
Using Claims~\ref{claim:mac1} and~\ref{claim-sis2}, we have $r(S_{d^*},V'_{d^*})\leq 2d^*$ and thus:

\begin{equation}\label{eq:star partial}
\forall v\in V'_{d^*}, d(v,C_{d^*})\leq 2d^*.
\end{equation}

Consider now a vertex $v\in V\setminus V'_{d^*}$ By definition of $V'_{d^*}$, it means that $d(v, \{a_i, i\in I_{d^*}^-\})\leq d^*$ and by definition of $ I_{d^*}^-$, it ensures $\exists i\in I_{d^*}^-, \forall u\in A_i, d(v,u)\leq 2d^*$.
This ensures:
\begin{equation}\label{eq:star complement}
\forall v\in V'_{d^*},  d(v,C_{d^*})\leq 2d^*.
\end{equation}

Equations~\ref{eq:star partial} and~\ref{eq:star complement} ensure $r(C_{d^*})\leq 2d^*$.
\qed \end{pf}

Claim~\ref{claim:mac3} and Equation~\ref{eq:tilde-star} imply $r(C_{\tilde d})\leq 2d^*$, which concludes the proof of Theorem~\ref{th: approx-mac_pcenter}.

\qed \end{pf}

We immediately deduce from Theorem~\ref{th: approx-mac_pcenter} and Proposition~\ref{prop:reduction} the main result of this section:

\begin{theoreme}\label{cor:approx}
For edge weighted graphs with lengths in $[\ell, 2\ell]$, \mar{Algorithm~\ref{algo: mac_pcenter} is a polynomial time approximation algorithm for {\tt Min P$p$CP} guaranteeing the ratio} $4\overline{deg}(G) +2$.
\end{theoreme}
\mar{In particular, on graphs with bounded degree, the ratio is constant:}

\mar{\begin{corollary}
{\tt Min P$p$CP} is constant approximable for graphs of bounded degree and edge lengths in $[\ell, 2\ell]$.
\end{corollary}}

\section{Conclusion}

In this paper, we strengthen the analysis of   {\tt Min P$p$CP} initiated in~\cite{DHM18}. In particular, in Section~\ref{sec:complexity}, we revisit the reduction we used in this previous paper to get a hardness result on planar graphs of bounded degree. The new reduction allows to prove that {\tt Min P$p$CP} is not approximable with a ratio less than $\frac{56}{55}$ on  subgrids of degree at most 3. Even thought the result does not generalize the one we previously obtained (the class is more restrictive but the new bound is closer to~1), the proof requires a much deeper analysis with techniques that might be useful for other problems. The main originality of our proof is the use of the intermediate graph  $\widetilde H_q$ (see Figure~\ref{fig:schema_transformation}): it can be seen as a perturbation of the subgrid $H_q$ that leads to a hard class for {\tt Min Vertex Cover}. 

Then, in Section~\ref{sec:approximation}, we propose some approximation results for this problem with, in particular,  a constant approximation for graphs of bounded degree and with edge lengths in $[\ell, 2\ell]$. To our knowledge, this is the first example of approximation for this problem and in addition it holds for a class of instances on which all our hardness results apply.  It provides a first gap between constant approximation ratios  and the hardness in approximation results we have obtained. Narrowing this gap for intermediate classes of graphs is a natural open question for further researches. 
In section~\ref{subsec: approx-pcenter}, we even show a stronger approximation result on trees.  However, we leave open the problem of whether {\tt Min P$p$CP} is NP-hard or polynomial on trees.

Most of our results apply for the uniform case only. Surprisingly,  Proposition~\ref{prop:reduction} and Theorem~\ref{cor:approx} are still valid for the case where edge lengths lie in $[\ell, 2\ell]$.  Finding polynomial cases and approximation results for {\tt Min P$p$CP} with general length system remains an important open question that would require new methods or tools.


Finally, when considering the feasibility conditions for {\tt Min P$p$CP}, we have introduced  the notion  of minimal articulation components (MACs) and the related {\tt Min MAC $p$-Center} Problem. We have shown that this problem is 2-approximable and that this is the best possible constant approximation ratio (Theorem~\ref{th: approx-mac_pcenter}). It is also polynomial on trees. Strengthening the 
\mar{study} of this notion and the complexity and approximation results for this problem on specific classes of instances is another question raised by the paper.

\def\bibname{Biblio}
\bibliographystyle{elsarticle-harv}
\bibliography{mybib}

\appendix
 
\refstepcounter{section}
\section*{Appendix}
\subsection{\bf Proof of Lemma~\ref{lem:relation_size_minVC}}

{\bf Lemma~\ref{lem:relation_size_minVC} }
\begin{em}
    Let $G=(V,E)$ be a graph and $G'=(V',E')$ be the graph obtained by inserting $2k_{uv}$ vertices on each edge $(u,v)\in E$, where $k_{uv}$ is a non-negative integer. Then we have
$$ \tau(G') = \tau(G) + \sum_{uv\in E} k_{uv}$$
\end{em}

\begin{pf}

For every edge $(u,v)\in E$ oriented from $u$ to $v$, denote $X_{uv}=\{x_{uv}^{1},\hdots, x_{uv}^{2k_{uv}}\}$ the set of vertices inserted on this edge. 
Note that at least $k_{uv}$ vertices are needed to cover vertices in $X_{uv}$.

Let $U \subset V$ a vertex cover of $G$: $\forall (u,v) \in E, \{u,v\} \cap U \neq \emptyset$.
We can build $U'\subset V'$ in $G'$ as follows. We initialize $U'$ with all vertices of $U$. Then, for every edge $(u,v)\in E$, if $u\in U$, we add vertices $x_{uv}^{2i}, 1\le i \le k_{uv}$ to $U'$.
Otherwise, $v\in U$ necessarily, then we add vertices $x_{uv}^{2i+1}, 0\le i \le k_{uv}-1$ to $U'$.
In both cases we have added exactly $k_{uv}$ vertices and all edges of $P^{G'}_{uv}$ are covered by $U'$, with $|U'|= |U| + \sum_{uv\in E} k_{uv}$.
Then $\tau(G')\le \tau(G) + k_{uv}$.

Assume now that $G'$ has a vertex cover $X'$.
For every $(u,v) \in E$, $P^{G'}_{uv}$ is covered by at least $k_{uv}+1$ vertices.
If $u,v\not\in X'$, we can transform $X'$ into $U'$ such that $u$ or $v$ is in $U'$.
Then $|U' \setminus V| \ge  \sum_{uv\in E} k_{uv}$.
Since at least one vertex between $u$ and $v$ is in $U = V \cap U'$, $U$ is a vertex cover for $G$.
Then $|U| = |U'| - k_{uv}$, thus $\tau(G)\le \tau(G') - k_{uv}$.

Hence $\tau(G')= \tau(G) + k_{uv}$ and the proof is complete.
\qed \end{pf}

\input{ListProblems}





\end{document}

%% file: ListProblems.tex
\subsection{\bf List of problems}
\label{app:listproblems}

\begin{problem} 
 \problemtitle{{\tt Min $p$-Center}}
 \probleminput{An edge-weighted graph $G=(V,E,L)$ and an integer $p$}
 \problemfeasible{Any $p$-center $C\subset V, |C|\leq p$}
 \problemoptimization{Minimize $r(C)=\max\limits_{v\in V}d(v,C)$. }
\end{problem}
\hskip 1cm 

\begin{problem} 
 \problemtitle{{\tt Min P$p$CP}}
 \probleminput{An edge-weighted graph $G=(V,E,L)$ and an integer $p$ ; the instance is denoted $(G,p)$}
 \problemfeasible{Any $p$-center $C\subset V, |C|\leq p$ satisfying ${\mathbb E}(C)<\infty$ }
 \problemoptimization{Minimize ${\mathbb E}(C)$. }
\end{problem}
\hskip 1cm

\begin{problem} 
 \problemtitle{{\tt Min MAC $p$-Center}}
 \probleminput{An edge-weighted graph $G=(V,E,L)$ and an integer $p$ }
 \problemfeasible{Any $p$-center $C\subset V, |C|\leq p$ satisfying ${\mathbb E}(C)<\infty$}
 \problemoptimization{Minimize $r(C)$. }
\end{problem}
\hskip 1cm 

\begin{problem} 
 \problemtitle{{\tt Min Partial $p$-Center}}
 \probleminput{An edge-weighted graph $G=(V,E,L)$, a subset $U\subset V$ and an integer $p$}
 \problemfeasible{Any $p$-center $C\subset U, |C|\leq p$}
 \problemoptimization{Minimize $r(C)$. }
\end{problem}

\hskip 1cm 

\begin{problem} 
 \problemtitle{{\tt Min Vertex Cover}}
 \probleminput{A graph $G=(V,E)$}
 \problemfeasible{A vertex cover i.e., a set $U\subseteq V$ such that every edge of $E$ is incident to at least one vertex of $U$}
 \problemoptimization{Minimize $|U|$. $\tau(G)$ denotes the minimum size of a vertex cover.}
\end{problem}

\hskip 1cm 

\begin{problem} 
 \problemtitle{{\tt Min Dominating Set}}
 \probleminput{A graph $G=(V,E)$}
 \problemfeasible{A dominating set i.e., a set $D\subseteq V$ such that every vertex of $V \setminus D$ is adjacent to a vertex of $D$.}
 \problemoptimization{Minimize $|D|$. $\gamma(G)$ denotes the minimum size of a dominating set.}
\end{problem}